# GAUSSIAN MODEL SELECTION WITH AN UNKNOWN VARIANCE

By Yannick Baraud, Christophe Giraud and Sylvie Huet

*Université de Nice Sophia Antipolis, Université de Nice Sophia Antipolis and INRA*

Let $Y$ be a Gaussian vector whose components are independent with a common unknown variance. We consider the problem of estimating the mean $\mu$ of $Y$ by model selection. More precisely, we start with a collection $\mathcal{S} = \{S_m, m \in \mathcal{M}\}$ of linear subspaces of $\mathbb{R}^n$ and associate to each of these the least-squares estimator of $\mu$ on $S_m$. Then, we use a data driven penalized criterion in order to select one estimator among these. Our first objective is to analyze the performance of estimators associated to classical criteria such as FPE, AIC, BIC and AMDL. Our second objective is to propose better penalties that are versatile enough to take into account both the complexity of the collection $\mathcal{S}$ and the sample size. Then we apply those to solve various statistical problems such as variable selection, change point detections and signal estimation among others. Our results are based on a nonasymptotic risk bound with respect to the Euclidean loss for the selected estimator. Some analogous results are also established for the Kullback loss.

**1. Introduction.** Let us consider the statistical model

$$(1.1) \qquad Y_i = \mu_i + \sigma \varepsilon_i, \qquad i = 1, \ldots, n,$$

where the parameters $\mu = (\mu_1, \ldots, \mu_n)' \in \mathbb{R}^n$ and $\sigma > 0$ are both unknown and the $\varepsilon_i$'s are i.i.d. standard Gaussian random variables. We want to estimate $\mu$ by model selection on the basis of the observation of $Y = (Y_1, \ldots, Y_n)'$.

To do this, we introduce a collection $\mathcal{S} = \{S_m, m \in \mathcal{M}\}$ of linear subspaces of $\mathbb{R}^n$, that hereafter will be called *models*, indexed by a finite or countable set $\mathcal{M}$. To each $m \in \mathcal{M}$ we can associate the least-squares estimator $\hat{\mu}_m = \Pi_m Y$ of $\mu$ relative to $S_m$ where $\Pi_m$ denotes the orthogonal projector onto $S_m$. Let us denote by $D_m$ the dimension of $S_m$ for $m \in \mathcal{M}$ and $\|\cdot\|$ the Euclidean norm on $\mathbb{R}^n$. The quadratic risk $\mathbb{E}[\|\mu - \hat{\mu}_m\|^2]$ of $\hat{\mu}_m$ with respect to this distance is given by

$$(1.2) \qquad \mathbb{E}[\|\mu - \hat{\mu}_m\|^2] = \inf_{s \in S_m} \|\mu - s\|^2 + D_m \sigma^2.$$

If we use this risk as a quality criterion, a best model is one minimizing the right-hand side of (1.2). Unfortunately, such a model is not available to the statistician since it depends on the unknown parameters $\mu$ and $\sigma^2$. A natural question then arises: to what extent can we select an element $\hat{m}(Y)$ of $\mathcal{M}$ depending on the data only, in such a way that the risk of the selected estimator $\hat{\mu}_{\hat{m}}$ be close to the minimal risk

$$(1.3) \qquad R(\mu, \mathcal{S}) = \inf_{m \in \mathcal{M}} \mathbb{E}[\|\mu - \hat{\mu}_m\|^2].$$

The art of *model selection* is to design such a selection rule in the best possible way. The standard way of solving the problem is to define $\hat{m}$ as the minimizer over $\mathcal{M}$ of some empirical criterion of the form

$$(1.4) \qquad \mathrm{Crit}_L(m) = \|Y - \Pi_m Y\|^2 \left(1 + \frac{\mathrm{pen}(m)}{n - D_m}\right)$$

or

$$(1.5) \qquad \mathrm{Crit}_K(m) = \frac{n}{2} \log\left(\frac{\|Y - \Pi_m Y\|^2}{n}\right) + \frac{1}{2} \mathrm{pen}'(m),$$

where pen and pen′ denote suitable (penalty) functions mapping $\mathcal{M}$ into $\mathbb{R}_+$. Note that these two criteria are equivalent (they select the same model) if pen and pen′ are related in the following way:

$$\mathrm{pen}'(m) = n \log\left(1 + \frac{\mathrm{pen}(m)}{n - D_m}\right), \quad \text{or} \quad \mathrm{pen}(m) = (n - D_m)(e^{\mathrm{pen}'(m)/n} - 1).$$

The present paper is devoted to investigating the performance of criterion (1.4) or (1.5) as a function of collection $\mathcal{S}$ and pen or pen′. More precisely, we want to deal with the following problems:

(P1) Given some collection $\mathcal{S}$ and an arbitrary nonnegative penalty function pen on $\mathcal{M}$, what will the performance $\mathbb{E}[\|\mu - \hat{\mu}_{\hat{m}}\|^2]$ of $\hat{\mu}_{\hat{m}}$ be?
(P2) What conditions on $\mathcal{S}$ and pen ensure that the ratio $\mathbb{E}[\|\mu - \hat{\mu}_{\hat{m}}\|^2]/R(\mu, \mathcal{S})$ is not too large.
(P3) Given a collection $\mathcal{S}$, what penalty should be recommended in view of minimizing (at least approximately) the risk of $\hat{\mu}_{\hat{m}}$?

It is beyond the scope of this paper to make an exhaustive historical review of the criteria of the form (1.4) and (1.5). We simply refer the interested reader to the first chapters of McQuarrie and Tsai (1998) for a nice and complete introduction to the domain. Let us only mention here some of the most popular criteria, namely FPE, AIC, BIC (or SIC) and AMDL which correspond respectively to the choices $\text{pen}(m) = 2D_m$, $\text{pen}'(m) = 2D_m$, $\text{pen}'(m) = D_m \log(n)$ and $\text{pen}'(m) = 3D_m \log(n)$. FPE was introduced in Akaike (1969) and is based on an unbiased estimate of the mean squared prediction error. AIC was proposed later by Akaike (1973) as a Kullback–Leibler information based model selection criterion. BIC and SIC are equivalent criteria which were respectively proposed by Schwarz (1978) and Akaike (1978) from a Bayesian perspective. More recently, Saito (1994) introduced AMDL as an information-theoretic based criterion. AMDL turns out to be a modified version of the Minimum Description Length criterion proposed by Rissanen (1983, 1984). The motivations for the construction of FPE, AIC, SIC and BIC criteria are a mixture of heuristic and asymptotic arguments. From both the theoretical and the practical point of view, these penalties suffer from the same drawback: their performance heavily depends on the sample size and the collection $\mathcal{S}$ at hand.

In recent years, more attention has been paid to the nonasymptotic point of view and a proper calibration of penalties taking into account the *complexity* (in a suitable sense) of the collection $\mathcal{S}$. A pioneering work based on the methodology of minimum complexity and dealing with discrete models and various stochastic frameworks including regression appeared in Barron and Cover (1991) and Barron (1991). It was then extended to various types of continuous models in Barron, Birgé and Massart (1999) and Birgé and Massart (1997, 2001a, 2007). Within the Gaussian regression framework, Birgé and Massart (2001a, 2007) consider model selection criteria of the form

$$(1.6) \qquad \text{crit}(m) = \|Y - \hat{\mu}_m\|^2 + \text{pen}(m)\sigma^2$$

and propose new penalty structures which depend on the complexity of the collection $\mathcal{S}$. These penalties can be viewed as generalizing Mallows' $C_p$ [heuristically introduced in Mallows (1973)] which corresponds to the choice $\text{pen}(m) = 2D_m$ in (1.6). However, Birgé and Massart only deal with the favorable situation where the variance $\sigma^2$ is known, although they provide some hints to estimate it in Birgé and Massart (2007).

Unlike Birgé and Massart, we consider here the more practical case where $\sigma^2$ is unknown. Yet our approach is similar in the sense that our objective is to propose new penalty structures for criteria (1.4) [or (1.5)] which allow us to take both the complexity of the collection and the sample size into account.

A possible application of the criteria we propose is variable selection in linear models. This problem has received a lot of attention in the literature. Recent development includes Tibshirani (1996) with the LASSO, Efron et al. (2004) with LARS, Candès and Tao (2007) for the Dantzig selector, Zou (2006) with the Adaptive LASSO, among others. Most of the recent literature assumes that $\sigma^2$ is known, or suitably estimated, and aim at designing an algorithm that solves the problem in polynomial time at the price of assumptions on the covariates to select. In contrast, our approach assumes nothing on $\sigma^2$ or the covariates, but requires that the number of these is not too large for a practical implementation.

The paper is organized as follows. In Section 2 we start with some examples of model selection problems among which variable selection, change point detection and denoising. This section gives the opportunity to both motivate our approach and make a review of some collections of models of interest. We address problem (P2) in Section 3 and analyze there FPE, AIC, BIC and AMDL criteria more specifically. In Section 4 we address problems (P1) and (P3) and introduce new penalty functions. In Section 5 we show how the statistician can take advantage of the flexibility of these new penalties to solve the model selection problems given in Section 2. Section 6 is devoted to two simulation studies allowing to assess the performances of our estimator. In the first one we consider the problem of detecting the nonzero components in the mean of a Gaussian vector and compare our estimator with BIC, AIC and AMDL. In the second study, we consider the variable selection problem and compare our procedure with the adaptive Lasso proposed by Zou (2006). In Section 7 we provide an analogue of our main result replacing the $\mathbb{L}^2$-loss by the Kullback loss. The remaining sections are devoted to the proofs.

To conclude this section, let us introduce some notation to be used throughout the paper. For each $m \in \mathcal{M}$, $D_m$ denotes the dimension of $S_m$, $N_m$ the quantity $n - D_m$ and $\mu_m = \Pi_m \mu$. We denote by $P_{\mu,\sigma^2}$ the distribution of $Y$. We endow $\mathbb{R}^n$ with the Euclidean inner product denoted $\langle \cdot, \cdot \rangle$. For all $x \in \mathbb{R}$, $(x)_+$ and $\lfloor x \rfloor$ denote respectively the positive and integer parts of $x$, and for $y \in \mathbb{R}$, $x \wedge y = \min\{x, y\}$ and $x \vee y = \max\{x, y\}$. Finally, we write $\mathbb{N}^*$ for the set of positive integers and $|m|$ for the cardinality of a set $m$.

**2. Some examples of model selection problems.** In order to illustrate and motivate the model selection approach to estimation, let us consider some examples of applications of practical interest. For each example, we shall describe the statistical problem at hand and the collection of models of interest. These collections will be characterized by a complexity index which is defined as follows.

DEFINITION 1. Let $M$ and $a$ be two nonnegative numbers. We say that a collection $\mathcal{S}$ of linear spaces $\{S_m, m \in \mathcal{M}\}$ has a finite complexity index $(M, a)$ if

$$|\{m \in \mathcal{M}, D_m = D\}| \le M e^{aD} \qquad \text{for all } D \ge 1.$$

Let us note here that not all countable families of models do have a finite complexity index.

2.1. *Detecting nonzero mean components.* The problem at hand is to recover the nonzero entries of a sparse high-dimensional vector $\mu$ observed with additional Gaussian noise. We assume that the vector $\mu$ in (1.1) has at most $p \leq n-2$ nonzero mean components but we do not know which are the null of these. Our goal is to find $m^* = \{i \in \{1, \ldots, n\} | \mu_i \neq 0\}$ and estimate $\mu$. Typically, $|m^*|$ is small as compared to the number of observations $n$. This problem has received a lot of attention in the recent years and various solutions have been proposed. Most of them rely on thresholding methods which require a suitable estimator of $\sigma^2$. We refer the interested reader to Abramovitch et al. (2006) and the references therein. Closer to our approach is the paper by Huet (2006) which is based on a penalized criterion related to AIC.

To handle this problem, we consider the set $\mathcal{M}$ of all subsets of $\{1, \ldots, n\}$ with cardinality not larger than $p$. For each $m \in \mathcal{M}$, we take for $S_m$ the linear space of those vectors $s$ in $\mathbb{R}^n$ such that $s_i = 0$ for $i \notin m$. By convention, $S_\varnothing = \{0\}$. Since the number of models with dimension $D$ is $\binom{n}{D} \leq n^D$, a complexity index for this collection is $(M, a) = (1, \log n)$.

2.2. *Variable selection.* Given a set of explanatory variables $x^{(1)}, \ldots, x^{(N)}$ and a response variable $y$ observed with additional Gaussian noise, we want to find a small subset of the explanatory variables that adequately explains $y$. This means that we observe $(Y_i, x_i^{(1)}, \ldots, x_i^{(N)})$ for $i = 1, \ldots, n$, where $x_i^{(j)}$ corresponds to the observation of the value of the variable $x^{(j)}$ in experiment number $i$, $Y_i$ is given by (1.1) and $\mu_i$ can be written as

$$\mu_i = \sum_{j=1}^{N} a_j x_i^{(j)},$$

where the $a_j$'s are unknown real numbers. Since we do not exclude the practical case where the number $N$ of explanatory variables is larger than the number $n$ of observations, this representation is not necessarily unique. We look for a subset $m$ of $\{1, \ldots, N\}$ such that the least-squares estimator $\hat{\mu}_m$ of $\mu$ based on the linear span $S_m$ of the vectors $x^{(j)} = (x_1^{(j)}, \ldots, x_n^{(j)})'$, $j \in m$, is as accurate as possible, restricting ourselves to sets $m$ of cardinality bounded by $p \leq n - 2$. By convention $S_\varnothing = \{0\}$.

A nonasymptotic treatment of this problem has been given by Birgé and Massart (2001a), Candès and Tao (2007) and Zou (2006) when $\sigma^2$ is known. To our knowledge, the practical case of an unknown value of $\sigma^2$ has not been analyzed from a nonasymptotic point of view. Note that when $N \geq n$ the traditional residual least-squares estimator cannot be used to estimate $\sigma^2$. Depending on our prior knowledge on the relative importance of the explanatory variables, we distinguish between two situations.

2.2.1. *A collection for "the ordered variable selection problem."* We consider here the favorable situation where the set of explanatory variables $\mathrm{x}^{(1)},\ldots,\mathrm{x}^{(p)}$ is ordered according to decreasing importance up to rank $p$ and introduce the collection

$$\mathcal{M}_o = \{\{1,\ldots,d\}, 1 \leq d \leq p\} \cup \varnothing,$$

subsets of $\{1,\ldots,N\}$. Since the collection contains at most one model per dimension, the family of models $\{S_m, m \in \mathcal{M}_o\}$ has a complexity index $(M,a) = (1,0)$.

2.2.2. *A collection for "the complete variable selection problem."* If we do not have much information about the relative importance of the explanatory variables $x^{(j)}$, it is more natural to choose for $\mathcal{M}$ the set of all subsets of $\{1,\ldots,N\}$ of cardinality not larger than $p$. For a given $D \geq 1$, the number of models with dimension $D$ is at most $\binom{N}{D} \leq N^D$ so that $(M,a) = (1, \log N)$ is a complexity index for the collection $\{S_m, m \in \mathcal{M}\}$.

2.3. *Change-points detection.* We consider the functional regression framework

$$Y_i = f(x_i) + \sigma\varepsilon_i, \qquad i = 1,\ldots,n,$$

where $\{x_1 = 0,\ldots,x_n\}$ is an increasing sequence of deterministic points of $[0,1)$ and $f$ an unknown real valued function on $[0,1)$. This leads to a particular instance of (1.1) with $\mu_i = f(x_i)$ for $i = 1,\ldots,n$. In such a situation, the loss function $\|\mu - \hat{\mu}\|^2 = \sum_{i=1}^n (f(x_i) - \hat{f}(x_i))^2$ is the discrete norm associated to the design $\{x_1,\ldots,x_n\}$.

We assume here that the unknown $f$ is either piecewise constant or piecewise linear with a number of change-points bounded by $p$. Our aim is to design an estimator $f$ which allows to estimate the number, locations and magnitudes of the jumps of either $f$ or $f'$, if any. The estimation of change-points of a function $f$ has been addressed by Lebarbier (2005) who proposed a model selection procedure related to Mallows' $C_p$.

2.3.1. *Models for detecting and estimating the jumps of $f$.* Since our loss function only involves the values of $f$ at the design points, natural models are those induced by piecewise constant functions with change-points among $\{x_2,\ldots,x_n\}$. A potential set $m$ of $q$ change-points is a subset $\{t_1,\ldots,t_q\}$ of $\{x_2,\ldots,x_n\}$ with $t_1 < \cdots < t_q$, $q \in \{0,\ldots,p\}$ with $p \leq n-3$, the set being empty when $q=0$. To a set $m$ of change-points $\{t_1,\ldots,t_q\}$ we associate the model
$$S_m = \{(g(x_1),\ldots,g(x_n))', g \in \mathcal{F}_m\},$$
where $\mathcal{F}_m$ is the space of piecewise constant functions of the form
$$\sum_{j=0}^{q} a_j \mathbb{1}_{[t_j,t_{j+1})} \qquad \text{with } (a_0,\ldots,a_q) \in \mathbb{R}^{q+1}, t_0 = x_1 \text{ and } t_{q+1} = 1,$$
so that the dimension of $S_m$ is $|m|+1$. Then we take for $\mathcal{M}$ the set of all subsets of $\{x_2,\ldots,x_n\}$ with cardinality bounded by $p$. For any $D$ with $1 \leq D \leq p+1$ the number of models with dimension $D$ is $\binom{n-1}{D-1} \leq n^D$ so that $(M,a) = (1,\log n)$ is a complexity index for this collection.

2.3.2. *A collection of models for detecting and estimating the jumps of $f'$.* Let us now turn to models for piecewise linear functions $g$ on $[0,1)$ with $q+1$ pieces so that $g'$ has at most $q \leq p$ jumps. We assume $p \leq n-4$. We denote by $\mathcal{C}([0,1))$ the set of continuous functions on $[0,1)$ and set $t_0 = 0$ and $t_{q+1} = 1$, as before. Given two nonnegative integers $j$ and $q$ such that $q <^j$, we set $\mathcal{D}_j = \{k2^{-j}, k=1,\ldots,2^j-1\}$ and define
$$\mathcal{M}_{j,q} = \{\{t_1,\ldots,t_q\} \subset \mathcal{D}_j, t_1 < \cdots < t_q\}$$
and
$$\mathcal{M} = \left(\bigcup_{j \geq 1} \bigcup_{q=1}^{(2^j-1)\wedge p} \mathcal{M}_{j,q}\right) \cup \{\varnothing\}.$$
For each $m = \{t_1,\ldots,t_q\} \in \mathcal{M}$ (with $m = \varnothing$ if $q=0$), we define $\mathcal{F}_m$ as the space of splines of degree 1 with knots in $m$, that is,
$$\mathcal{F}_m = \left\{\sum_{k=0}^{q}(\alpha_k x + \beta_k)\mathbb{1}_{[t_k,t_{k+1})}(x) \in \mathcal{C}([0,1)), (\alpha_k,\beta_k)_{0 \leq k \leq q} \in \mathbb{R}^{2(q+1)}\right\}$$
and the corresponding model
$$S_m = \{(g(x_1),\ldots,g(x_n))', g \in \mathcal{F}_m\} \subset \mathbb{R}^n.$$
Note that $2 \leq \dim(S_m) \leq \dim(\mathcal{F}_m) = |m|+2$ because of the continuity constraint. Besides, let us observe that $\mathcal{M}$ is countable and that the number of models $S_m$ with a dimension $D$ in $\{1,\ldots,p+2\}$ is infinite. This implies that the collection has no (finite) complexity index.

2.4. *Estimating an unknown signal.* We consider the problem of estimating a (possibly) anisotropic signal in $\mathbb{R}^d$ observed at discrete times with additional noise. This means that we observe the vector $Y$ given by (1.1) with

$$\mu_i = f(x_i), \qquad i = 1, \ldots, n, \tag{2.1}$$

where $x_1, \ldots, x_n \in [0,1)^d$ and $f$ is an unknown function mapping $[0,1)^d$ into $\mathbb{R}$. To estimate $f$ we use models of piecewise polynomial functions on partitions of $[0,1)^d$ into hyperrectangles. We consider the set of indices

$$\mathcal{M} = \{(r, k_1, \ldots, k_d), r \in \mathbb{N}, k_1, \ldots, k_d \in \mathbb{N}^* \text{ with } (r+1)^d k_1 \cdots k_d \leq n - 2\}.$$

For $m = (r, k_1, \ldots, k_d) \in \mathcal{M}$, we set $\mathcal{J}_m = \prod_{i=1}^d \{1, \ldots, k_i\}$ and denote by $\mathcal{F}_m$ the space of piecewise polynomials $P$ such that the restriction of $P$ to each hyperrectangle $\prod_{i=1}^d [(j_i - 1)k_i^{-1}, j_i k_i^{-1})$ with $j \in \mathcal{J}_m$ is a polynomial in $d$ variables of degree not larger than $r$. Finally, we consider the collection of models

$$S_m = \{(P(x_1), \ldots, P(x_n))', P \in \mathcal{F}_m\}, \qquad m \in \mathcal{M}.$$

Note that when $m = (r, k_1, \ldots, k_d)$, the dimension of $S_m$ is not larger than $(r+1)^d k_1 \cdots k_d$. A similar collection of models was introduced in Barron, Birgé and Massart (1999) for the purpose of estimating a density on $[0,1)^d$ under some Hölderian assumptions.

**3. Analyzing penalized criteria with regard to family complexity.** Throughout the section, we set $\phi(x) = (x - 1 - \log(x))/2$ for $x \geq 1$ and denote by $\phi^{-1}$ the reciprocal of $\phi$. We assume that the collection of models satisfies for some $K > 1$ and $(M, a) \in \mathbb{R}_+^2$ the following assumption.

ASSUMPTION ($H_{K,M,a}$). *The collection of models $\mathcal{S} = \{S_m, m \in \mathcal{M}\}$ has a complexity index $(M, a)$ and satisfies*

$$\forall m \in \mathcal{M}, D_m \leq D_{\max} = \lfloor (n - \gamma_1)_+ \rfloor \wedge \lfloor ((n+2)\gamma_2 - 1)_+ \rfloor,$$

*where*

$$\gamma_1 = (2t_{a,K}) \vee \frac{t_{a,K} + 1}{t_{a,K} - 1},$$

$$\gamma_2 = \frac{2\phi(K)}{(t_{a,K} - 1)^2}$$

*and*

$$t_{a,K} = K\phi^{-1}(a) > 1.$$

If $a = 0$ and $a = \log(n)$, Assumption $(\mathrm{H}_{K,M,a})$ amounts to assuming $D_m \leq \delta(K)n$ and $D_m \leq \delta(K)n/\log^2(n)$, respectively, for all $m \in \mathcal{M}$ where $\delta(K) < 1$ is some constant depending on $K$ only. In any case, since $\gamma_2 \leq 2\phi(K)(K-1)^{-2} \leq 1/2$, Assumption $(\mathrm{H}_{K,M,a})$ implies that $D_{\max} \leq n/2$.

3.1. *Bounding the risk of $\hat{\mu}_{\hat{m}}$ under penalty constraints.* The following holds.

THEOREM 1. *Let $K > 1$ and $(M, a) \in \mathbb{R}_+^2$. Assume that the collection $\mathcal{S} = \{S_m, m \in \mathcal{M}\}$ satisfies $(\mathrm{H}_{K,M,a})$. If $\hat{m}$ is selected as a minimizer of $\mathrm{Crit}_L$ [defined by (1.4)] among $\mathcal{M}$ and if* pen *satisfies*

$$\text{(3.1)} \qquad \mathrm{pen}(m) \geq K^2 \phi^{-1}(a) D_m \qquad \forall m \in \mathcal{M},$$

*then the estimator $\hat{\mu}_{\hat{m}}$ satisfies*

$$\text{(3.2)} \quad \mathbb{E}\left[\frac{\|\mu - \hat{\mu}_{\hat{m}}\|^2}{\sigma^2}\right]$$
$$\leq \frac{K}{K-1} \inf_{m \in \mathcal{M}} \left[\frac{\|\mu - \mu_m\|^2}{\sigma^2}\left(1 + \frac{\mathrm{pen}(m)}{n - D_m}\right) + \mathrm{pen}(m) - D_m\right] + R,$$

*where*

$$R = \frac{K}{K-1}\left[K^2 \phi^{-1}(a) + 2K + \frac{8KMe^{-a}}{(e^{\phi(K)/2} - 1)^2}\right].$$

*In particular, if $\mathrm{pen}(m) = K^2 \phi^{-1}(a) D_m$ for all $m \in \mathcal{M}$,*

$$\text{(3.3)} \qquad \mathbb{E}[\|\mu - \hat{\mu}_{\hat{m}}\|^2] \leq C\phi^{-1}(a)[R(\mu, \mathcal{S}) \vee \sigma^2]$$

*where $C$ is a constant depending on $K$ and $M$ only and $R(\mu, \mathcal{S})$ the quantity defined at equation (1.3).*

If we exclude the situation where $\{0\} \in \mathcal{S}$, one has $R(\mu, \mathcal{S}) \geq \sigma^2$. Then, (3.3) shows that the choice $\mathrm{pen}(m) = K^2 \phi^{-1}(a) D_m$ leads to a control of the ratio $\mathbb{E}[\|\mu - \hat{\mu}_{\hat{m}}\|^2]/R(\mu, \mathcal{S})$ by the quantity $C\phi^{-1}(a)$ which only depends on $K$ and the complexity index $(M, a)$. For a typical collection of models, $a$ is either of order of a constant (independent of $n$) or of order of a $\log(n)$. In the first case, the risk bound we get leads to an oracle-type inequality showing that the resulting estimator achieves up to constant the best trade-off between the bias and the variance term. In the second case, $\phi^{-1}(a)$ is of order of a $\log(n)$ and the risk of the estimator differs from $R(\mu, \mathcal{S})$ by a logarithmic factor. For the problem described in Section 2.1, this extra logarithmic factor is known to be unavoidable [see Donoho and Johnstone (1994), Theorem 3]. We shall see in Section 3.3 that the constraint (3.1) is sharp at least in the typical situations where $a = 0$ and $a = \log(n)$.

3.2. *Analysis of some classical penalties with regard to complexity.* In the sequel, we make a review of classical penalties and analyze their performance in the light of Theorem 1.

*FPE and AIC.* As already mentioned, FPE corresponds to the choice $\text{pen}(m) = 2D_m$. If the complexity index $a$ belongs to $[0, \phi(2))$ $[\phi(2) \approx 0.15]$, then this penalty satisfies (3.1) with $K = \sqrt{2/\phi^{-1}(a)} > 1$. If the complexity index of the collection is $(M, a) = (1, 0)$, by assuming that

$$D_m \leq \min\{n - 6, 0.39(n+2) - 1\}$$

we ensure that Assumption $(\text{H}_{K,M,a})$ holds and we deduce from Theorem 1 that (3.2) is satisfied with $K/(K-1) < 3.42$. For such collections, the use of FPE leads thus to an oracle-type inequality. The AIC criterion corresponds to the penalty $\text{pen}(m) = N_m(e^{2D_m/n} - 1) \geq 2N_m D_m/n$ and has thus similar properties provided that $N_m/n$ remains bounded from below by some constant larger than $1/2$.

*AMDL and BIC.* The AMDL criterion corresponds to the penalty

$$(3.4) \qquad \text{pen}(m) = N_m(e^{3D_m \log(n)/n} - 1) \geq 3N_m n^{-1} D_m \log(n).$$

This penalty can cope with the (complex) collection of models introduced in Section 2.1 for the problem of detecting the nonzero mean components in a Gaussian vector. In this case, the complexity index of the collection can be taken as $(M, a) = (1, \log(n))$ and since $\phi^{-1}(a) \leq 2\log(n)$, inequality (3.1) holds with $K = \sqrt{2}$. As soon as for all $m \in \mathcal{M}$,

$$(3.5) \qquad D_m \leq \min\left\{n - 5.7\log(n), \frac{0.06(n+2)}{(3\log(n)-1)^2} - 1\right\},$$

Assumption $(\text{H}_{K,M,a})$ is fulfilled and $\hat{\mu}_{\hat{m}}$ then satisfies (3.2) with $K/(K-1) < 3.42$. Actually, this result has an asymptotic flavor since (3.5) and therefore $(\text{H}_{K,M,a})$ hold for very large values of $n$ only. For a more practical point of view, we shall see in Section 6 that AMDL penalty is too large and thus favors small dimensional linear spaces too much. The BIC criterion corresponds to the choice $\text{pen}(m) = N_m(e^{D_m \log(n)/n} - 1)$ and one can check that $\text{pen}(m)$ stays smaller than $\phi^{-1}(\log(n))D_m$ when $n$ is large. Consequently, Theorem 1 cannot justify the use of the BIC criterion for the collection above. In fact, we shall see in the next section that BIC is inappropriate in this case.

When the complexity parameter $a$ is independent of $n$, criteria AMDL and BIC satisfy (3.1) for $n$ large enough. Nevertheless, the logarithmic factor involved in these criteria has the drawback to overpenalize large dimensional linear spaces. One consequence is that the risk bound (3.2) differs from an oracle inequality by a logarithmic factor.

3.3. *Minimal penalties.* The aim of this section is to show that the constraint (3.1) on the size of the penalty is sharp. We shall restrict ourselves to the cases where $a = 0$ and $a = \log(n)$. Similar results have been established in Birgé and Massart (2007) for criteria of the form (1.6). The interested reader can find the proofs of the following propositions in Baraud, Giraud and Huet (2007).

3.3.1. *Case $a = 0$.* For collections with such a complexity index, we have seen that the conditions of Theorem 1 are fulfilled as soon as $\text{pen}(m) \geq CD_m$ for all $m$ and some universal constant $C > 1$. Besides, the choice of penalties of the form $\text{pen}(m) = CD_m$ for all $m$ leads to oracle inequalities. The following proposition shows that the constraint $C > 1$ is necessary to avoid the overfitting phenomenon.

PROPOSITION 1. *Let $\mathcal{S} = \{S_m, m \in \mathcal{M}\}$ be a collection of models with complexity index $(1, 0)$. Assume that $\text{pen}(\bar{m}) < D_{\bar{m}}$ for some $\bar{m} \in \mathcal{M}$ and set $C = \text{pen}(\bar{m})/D_{\bar{m}}$. If $\mu = 0$, the index $\hat{m}$ which minimizes criterion (1.4) satisfies*

$$\mathbb{P}\left(D_{\hat{m}} \geq \frac{1-C}{2} D_{\bar{m}}\right) \geq 1 - ce^{-c'(N_{\bar{m}} \wedge D_{\bar{m}})},$$

*where $c$ and $c'$ are positive functions of $C$ only.*

Explicit values of $c$ and $c'$ can be found in the proof.

3.3.2. *Case $a = \log(n)$.* We restrict ourselves to the collection described in Section 2.1. We have already seen that the choice of penalties of the form $\text{pen}(m) = 2CD_m \log n$ for all $m$ with $C > 1$ was leading to a nearly optimal bias and variance trade-off [up to an unavoidable $\log(n)$ factor] in the risk bounds. We shall now see that the constraint $C > 1$ is sharp.

PROPOSITION 2. *Let $C_0 \in ]0, 1[$. Consider the collection of linear spaces $\mathcal{S} = \{S_m | m \in \mathcal{M}\}$ described in Section 2.1, and assume that $p \leq (1 - C_0)n$ and $n > e^{2/C_0}$. Let $\text{pen}$ be a penalty satisfying $\text{pen}(m) \leq 2C_0^4 D_m \log(n)$ for all $m \in \mathcal{M}$. If $\mu = 0$, the cardinality of the subset $\hat{m}$ selected as a minimizer of criterion (1.4) satisfies*

$$\mathbb{P}(|\hat{m}| > \lfloor(1-C_0)D\rfloor) \geq 1 - 2\exp\left(-c\frac{n^{1-C_0}}{\sqrt{\log(n)}}\right),$$

*where $D = \lfloor c'n^{1-C_0}/\log^{3/2}(n)\rfloor \wedge p$ and $c$, $c'$ are positive functions of $C_0$ (to be explicitly given in the proof).*

Proposition 2 shows that AIC and FPE should not be used for model selection purposes with the collection of Section 2.1. Moreover, if $p\log(n)/n \leq \kappa < \log(2)$ then the BIC criterion satisfies

$$\mathrm{pen}(m) = N_m(e^{D_m \log(n)/n} - 1) \leq e^\kappa D_m \log(n) < 2 D_m \log(n)$$

and also appears inadequate to cope with the complexity of this collection.

**4. From general risk bounds to new penalized criteria.** Given an arbitrary penalty pen, our aim is to establish a risk bound for the estimator $\hat{\mu}_{\hat{m}}$ obtained from the minimization of (1.4). The analysis of this bound will lead us to propose new penalty structures that take into account the complexity of the collection. Throughout this section we shall assume that $D_m \leq n-2$ for all $m \in \mathcal{M}$.

The main theorem of this section uses the function Dkhi defined below.

DEFINITION 2. Let $D, N$ be two positive numbers and $X_D, X_N$ be two independent $\chi^2$ random variables with degrees of freedom $D$ and $N$ respectively. For $x \geq 0$, we define

$$(4.1) \qquad \mathsf{Dkhi}[D, N, x] = \frac{1}{\mathbb{E}(X_D)} \times \mathbb{E}\left[\left(X_D - x\frac{X_N}{N}\right)_+\right].$$

Note that for $D$ and $N$ fixed, $x \mapsto \mathsf{Dkhi}[D, N, x]$ is decreasing from $[0, +\infty)$ into $(0, 1]$ and satisfies $\mathsf{Dkhi}[D, N, 0] = 1$.

THEOREM 2. *Let $\mathcal{S} = \{S_m, m \in \mathcal{M}\}$ be some collection of models such that $N_m \geq 2$ for all $m \in \mathcal{M}$. Let pen be an arbitrary penalty function mapping $\mathcal{M}$ into $\mathbb{R}^+$. Assume that there exists an index $\hat{m}$ among $\mathcal{M}$ which minimizes (1.4) with probability 1. Then, the estimator $\hat{\mu}_{\hat{m}}$ satisfies for all constants $c \geq 0$ and $K > 1$,*

$$\mathbb{E}\left[\frac{\|\mu - \hat{\mu}_{\hat{m}}\|^2}{\sigma^2}\right]$$

$$(4.2) \qquad \leq \frac{K}{K-1} \inf_{m \in \mathcal{M}} \left[\frac{\|\mu - \mu_m\|^2}{\sigma^2}\left(1 + \frac{\mathrm{pen}(m)}{N_m}\right) + \mathrm{pen}(m) - D_m\right]$$

$$+ \Sigma,$$

*where*

$$\Sigma = \frac{Kc}{K-1}$$

$$+ \frac{2K^2}{K-1} \sum_{m \in \mathcal{M}} (D_m + 1)\mathsf{Dkhi}\left[D_m + 1, N_m - 1, \frac{N_m - 1}{KN_m}(\mathrm{pen}(m) + c)\right].$$

Note that a minimizer of $\operatorname{Crit}_L$ does not necessarily exist for an arbitrary penalty function, unless $\mathcal{M}$ is finite. Take for example, $\mathcal{M} = \mathbb{Q}^n$ and for all $m \in \mathcal{M}$ set $\operatorname{pen}(m) = 0$ and $S_m$ the linear span of $m$. Since $\inf_{m \in \mathcal{M}} \|Y - \Pi_m Y\|^2 = 0$ and $Y \notin \bigcup_{m \in \mathcal{M}} S_m$ a.s., $\hat{m}$ does not exist with probability 1. In the case where $\hat{m}$ does exist with probability 1, the quantity $\Sigma$ appearing in right-hand side of (4.2) can either be calculated numerically or bounded by using Lemma 6 below.

Let us now turn to an analysis of inequality (4.2). Note that the right-hand side of (4.2) consists of the sum of two terms,

$$\frac{K}{K-1} \inf_{m \in \mathcal{M}} \left[ \frac{\|\mu - \mu_m\|^2}{\sigma^2} \left(1 + \frac{\operatorname{pen}(m)}{N_m}\right) + \operatorname{pen}(m) - D_m \right]$$

and $\Sigma = \Sigma(\operatorname{pen})$, which vary in opposite directions with the size of pen. There is clearly no hope in optimizing this sum with respect to pen without any prior information on $\mu$. Since only $\Sigma$ depends on known quantities, we suggest choosing the penalty in view of controlling its size. As already seen, the choice $\operatorname{pen}(m) = K^2 \phi^{-1}(a) D_m$ for some $K > 1$ allows us to obtain a control of $\Sigma$ which is independent of $n$. This choice has the following drawbacks. First, the penalty penalizes the same all the models of a given dimension, although one could wish to associate a smaller penalty to some of these because they possess a simpler structure. Second, it turns out that in practice these penalties are a bit too large and leads to an underfitting of the true by advantaging too much small dimensional models. In order to avoid these drawbacks, we suggest to use the penalty structures introduced in the next section.

4.1. *Introducing new penalty functions.* We associate to the collection of models $\mathcal{S}$ a collection $\mathcal{L} = \{L_m, m \in \mathcal{M}\}$ of nonnegative numbers (weights) such that

$$(4.3) \qquad \Sigma' = \sum_{m \in \mathcal{M}} (D_m + 1) e^{-L_m} < +\infty.$$

When $\Sigma' = 1$ then the choice of sequence $\mathcal{L}$ can be interpreted as a choice of a prior distribution $\pi$ on the set $\mathcal{M}$. This a priori choice of a collection of $L_m$'s gives a Bayesian flavor to the selection rule. We shall see in the next section how the sequence $\mathcal{L}$ can be chosen in practice according to the collection at hand.

DEFINITION 3. For $0 < q \leq 1$ we define $\mathsf{EDkhi}[D, N, q]$ as the unique solution of the equation $\mathsf{Dkhi}[D, N, \mathsf{EDkhi}[D, N, q]] = q$.

Given some $K > 1$, let us define the penalty function $\operatorname{pen}_{K,\mathcal{L}}$

$$(4.4) \quad \operatorname{pen}_{K,\mathcal{L}}(m) = K \frac{N_m}{N_m - 1} \mathsf{EDkhi}[D_m + 1, N_m - 1, e^{-L_m}] \qquad \forall m \in \mathcal{M}.$$

PROPOSITION 3. *If* $\text{pen} = \text{pen}_{K,\mathcal{L}}$ *for some sequence of weights $\mathcal{L}$ satisfying (4.3), then there exists an index $\hat{m}$ among $\mathcal{M}$ which minimizes (1.4) with probability 1. Besides, the estimator $\hat{\mu}_{\hat{m}}$ satisfies (4.2) with $\Sigma \leq 2K^2 \Sigma'/(K-1)$.*

As we shall see in Section 6.1, the penalty $\text{pen}_{K,\mathcal{L}}$ or at least an upper bound can easily be computed in practice. From a more theoretical point of view, an upper bound for $\text{pen}_{K,\mathcal{L}}(m)$ is given in the following proposition, the proof of which is postponed to Section 10.2.

PROPOSITION 4. *Let $m \in \mathcal{M}$ such that $N_m \geq 7$ and $D_m \geq 1$. We set $D = D_m + 1$, $N = N_m - 1$ and*

$$\Delta = \frac{L_m + \log 5 + 1/N}{1 - 5/N}.$$

*Then, we have the following upper bound on the penalty $\text{pen}_{K,\mathcal{L}}(m)$:*

$$(4.5) \quad \text{pen}_{K,\mathcal{L}}(m) \leq \frac{K(N+1)}{N}\left[1 + e^{2\Delta/(N+2)}\sqrt{\left(1 + \frac{2D}{N+2}\right)\frac{2\Delta}{D}}\right]^2 D.$$

*When $D_m = 0$ and $N_m \geq 4$, we have the upper bound*

$$(4.6) \quad \text{pen}_{K,\mathcal{L}}(m) \leq \frac{3K(N+1)}{N}\left[1 + e^{2L_m/N}\sqrt{\left(1 + \frac{6}{N}\right)\frac{2L_m}{3}}\right]^2.$$

*In particular, if $L_m \vee D_m \leq \kappa n$ for some $\kappa < 1$, then there exists a constant $C$ depending on $\kappa$ and $K$ only, such that*

$$\text{pen}_{K,\mathcal{L}}(m) \leq C(L_m \vee D_m)$$

*for any $m \in \mathcal{M}$.*

We derive from Proposition 4 and Theorem 2 (with $c = 0$) the following risk bound for the estimator $\hat{\mu}_{\hat{m}}$.

COROLLARY 1. *Let $\kappa < 1$. If for all $m \in \mathcal{M}$, $N_m \geq 7$ and $L_m \vee D_m \leq \kappa n$, then $\hat{\mu}_{\hat{m}}$ satisfies*

$$(4.7) \quad \mathbb{E}\left[\frac{\|\mu - \hat{\mu}_{\hat{m}}\|^2}{\sigma^2}\right] \leq C\left[\inf_{m \in \mathcal{M}}\left\{\frac{\|\mu - \mu_m\|^2}{\sigma^2} + D_m \vee L_m\right\} + \Sigma'\right],$$

*where $C$ is a positive quantity depending on $\kappa$ and $K$ only.*

Note that (4.7) turns out to be an oracle-type inequality as soon as one can choose $L_m$ of the order of $D_m$ for all $m$. Unfortunately, this is not always possible if one wants to keep the size of $\Sigma'$ under control. Finally, let us mention that the structure of our penalties, $\text{pen}_{K,\mathcal{L}}$, is flexible enough to recover any penalty function pen by choosing the family of weights $\mathcal{L}$ adequately. Namely, it suffices to take

$$L_m = -\log\left(\mathsf{Dkhi}\left[D_m+1, N_m-1, \frac{(N_m-1)\,\text{pen}(m)}{KN_m}\right]\right)$$

to obtain $\text{pen}_{K,\mathcal{L}} = \text{pen}$. Nevertheless, this choice of $\mathcal{L}$ does not ensure that (4.3) holds true unless $\mathcal{M}$ is finite.

## 5. How to choose the weights.

5.1. *One simple way.* One can proceed as follows. If the complexity index of the collection is given by the pair $(M, a)$, then the choice

(5.1) $$L_m = a' D_m \qquad \forall m \in \mathcal{M}$$

for some $a' > a$ leads to the following control of $\Sigma'$:

$$\Sigma' \le M \sum_{D \ge 1} D e^{-(a'-a)(D-1)} = M(1 - e^{-(a'-a)})^{-2}.$$

In practice, this choice of $\mathcal{L}$ is often too rough. One of its nonattractive features lies in the fact that the resulting penalty penalizes the same all the models of a given dimension. Since it is not possible to give a universal recipe for choosing the sequence $\mathcal{L}$, in the sequel we consider the examples presented in Section 2 and in each case motivate a choice of a specific sequence $\mathcal{L}$ by theoretical or practical considerations.

5.2. *Detecting nonzero mean components.* For any $D \in \{0, \ldots, p\}$ and $m \in \mathcal{M}$ such that $|m| = D$, we set

$$L_m = L(D) = \log\left[\binom{n}{D}\right] + 2\log(D+1)$$

and $\text{pen}(m) = \text{pen}_{K,\mathcal{L}}(m)$ where $K$ is some fixed constant larger than 1. Since $\text{pen}(m)$ only depends on $|m|$, we write

(5.2) $$\text{pen}(m) = \overline{\text{pen}}(|m|).$$

From a practical point of view, $\hat{m}$ can be computed as follows. Let $Y_{(n)}^2, \ldots, Y_{(1)}^2$ be random variables obtained by ordering $Y_1^2, \ldots, Y_n^2$ in the following way:

$$Y_{(n)}^2 < Y_{(n-1)}^2 < \cdots < Y_{(1)}^2 \qquad \text{a.s.}$$

and $\hat{D}$ the integer minimizing over $D \in \{0, \ldots, p\}$ the quantity

$$\sum_{i=D+1}^{n} Y_{(i)}^2 \left(1 + \frac{\overline{\text{pen}}(D)}{n-D}\right). \tag{5.3}$$

Then the subset $\hat{m}$ coincides with $\{(1), \ldots, (\hat{D})\}$ if $\hat{D} \geq 1$ and $\varnothing$ otherwise. In Section 6 a simulation study evaluates the performance of this method for several values of $K$.

From a theoretical point of view, our choice of $L_m$'s implies the following bound on $\Sigma'$:

$$\Sigma' \leq \sum_{D=0}^{p} \binom{n}{D}(D+1)e^{-L(D)}$$

$$\leq \sum_{D=1}^{p} \frac{1}{D} \leq 1 + \log(p+1) \leq 1 + \log(n).$$

As to the penalty, let us fix some $m$ in $\mathcal{M}$ with $|m| = D$. The usual bound $\log[\binom{n}{D}] \leq D \log(n)$ implies $L_m \leq D(2 + \log n) \leq p(2 + \log(n))$ and consequently, under the assumption

$$p \leq \frac{\kappa n}{2 + \log n} \wedge (n-7)$$

for some $\kappa < 1$, we deduce from Corollary 1 that for some constant $C' = C'(\kappa, K)$, the estimator $\hat{\mu}_{\hat{m}}$ satisfies

$$\mathbb{E}[\|\mu - \hat{\mu}_{\hat{m}}\|^2] \leq C' \inf_{m \in \mathcal{M}} [\|\mu - \mu_m\|^2 + (D_m + 1)\log(n)\sigma^2]$$

$$\leq C'(1 + |m^*|)\log(n)\sigma^2.$$

As already mentioned, we know that the $\log(n)$ factor in the risk bound is unavoidable. Unlike the former choice of $\mathcal{L}$ suggested by (5.1) [with $a' = \log(n) + 1$, e.g.], the bound for $\Sigma'$ we get here is not independent of $n$ but rather grows with $n$ at rate $\log(n)$. As compared to the former, this latter weighting strategy leads to similar risk bounds and to a better performance of the estimator in practice.

5.3. *Variable selection.* We propose to handle simultaneously complete and ordered variable selection. First, we consider the $p$ explanatory variables that we believe to be the most important among the set of the $N$ possible ones. Then, we index these from 1 to $p$ by decreasing order of importance and index those $N-p$ remaining ones arbitrarily. We do not assume that our guess on the importance of the various variables is right or not. We define $\mathcal{M}_o$ and $\mathcal{M}$ according to Section 2.2 and for some $c > 0$ set $L_m = c|m|$, if $m \in \mathcal{M}_o$, and otherwise set

$$L_m = L(|m|) \qquad \text{where } L(D) = \log\left[\binom{N}{D}\right] + \log p + \log(D+1).$$

For $K > 1$, we select the subset $\hat{m}$ as the minimizer among $\mathcal{M}$ of the criterion $m \mapsto \text{Crit}_L(m)$ given by (1.4) with $\text{pen}(m) = \text{pen}_{K,\mathcal{L}}(m)$. Except in the favorable situation where the vectors $x^{(j)}$ are orthogonal in $\mathbb{R}^n$ there seems, unfortunately, to be no way of computing $\hat{m}$ in polynomial time. Nevertheless, the method can be applied for reasonable values of $N$ and $p$ as shown in Section 6.3. From a theoretical point of view, our choice of $L_m$'s leads to the following bound on the residual term $\Sigma'$:

$$\Sigma' \leq \sum_{m \in \mathcal{M}_o} (|m|+1)e^{-L_m} + \sum_{m \in \mathcal{M} \setminus \mathcal{M}_o} (|m|+1)e^{-L_m}$$

$$\leq \sum_{D=0}^{p} (D+1)e^{-cD} + \sum_{D=1}^{p} \binom{D}{N}(D+1)e^{-L(D)}$$

$$\leq 1 + (1-e^{-c})^{-2}.$$

Besides, we deduce from Corollary 1 that if $p$ satisfies

$$p \leq \frac{\kappa n}{c} \wedge \frac{\kappa n}{2 + \log N} \wedge (n-7) \qquad \text{with } \kappa < 1,$$

then

(5.4) $$\mathbb{E}[\|\mu - \hat{\mu}_{\hat{m}}\|^2] \leq C(\kappa, K, c)(B_o \wedge B_c),$$

where

$$B_o = \inf_{m \in \mathcal{M}_o}(\|\mu - \mu_m\|^2 + (|m|+1)\sigma^2),$$

$$B_c = \inf_{m \in \mathcal{M}}[\|\mu - \mu_m\|^2 + (|m|+1)\log(eN)\sigma^2].$$

It is interesting to compare the risk bound (5.4) with the one we can get by using the former choice of weights $\mathcal{L}$ given in (5.1) [with $a' = \log(N) + 1$], that is

(5.5) $$\mathbb{E}[\|\mu - \hat{\mu}_{\hat{m}}\|^2] \leq C'(\kappa, K) B_c.$$

Up to constants, we see that (5.4) improves (5.5) by a $\log(N)$ factor whenever the minimizer $m^*$ of $\mathbb{E}[\|\mu - \hat{\mu}_m\|^2]$ among $\mathcal{M}$ does belong to $\mathcal{M}_o$.

5.4. *Multiple change-points detection.* In this section, we consider the problems of change-points detection presented in Section 2.3.

5.4.1. *Detecting and estimating the jumps of $f$.* We consider here the collection of models described in Section 2.3.1 and associate to each $m$ the weight $L_m$ given by

$$L_m = L(|m|) = \log\left[\binom{n-1}{|m|}\right] + 2\log(|m|+2),$$

where $K$ is some number larger than 1. This choice gives the following control on $\Sigma'$:

$$\Sigma' = \sum_{D=0}^{p}\binom{n-1}{D}(D+2)e^{-L(D)} = \sum_{D=0}^{p}\frac{1}{D+2} \leq \log(p+2).$$

Let $D$ be some arbitrary positive integer not larger than $p$. If $f$ belongs to the class of functions which are piecewise constant on an arbitrary partition of $[0,1)$ into $D$ intervals, then $\mu = (f(x_1),\ldots,f(x_n))'$ belongs to some $S_m$ with $m \in \mathcal{M}$ and $|m| \leq D-1$. We deduce from Corollary 1 that if $p$ satisfies

$$p \leq \frac{\kappa n - 2}{2 + \log n} \wedge (n-8)$$

for some $\kappa < 1$, then

$$\mathbb{E}[\|\mu - \hat{\mu}_{\hat{m}}\|^2] \leq C(\kappa, K) D \log(n) \sigma^2.$$

5.4.2. *Detecting and estimating the jumps of $f'$.* In this section, we deal with the collection of models of Section 2.3.2. Note that this collection is not finite. We use the following weighting strategy. For any pair of integers $j, q$ such that $q \leq 2^j - 1$, we set

$$L(j, q) = \log\left[\binom{2^j - 1}{q}\right] + q + 2\log j.$$

Since an element $m \in \mathcal{M}$ may belong to different $\mathcal{M}_{j,q}$, we set $L_m = \inf\{L(j, q), m \in \mathcal{M}_{j,q}\}$. This leads to the following control of $\Sigma'$:

$$\Sigma' \leq \sum_{j \geq 1} \sum_{q=0}^{(2^j-1)\wedge p} |\mathcal{M}_{j,q}|(q+3)\frac{e^{-q}}{\binom{2^j-1}{q}j^2}$$

$$\leq \sum_{j \geq 1} \frac{1}{j^2} \sum_{q \geq 0}(q+3)e^{-q}$$

$$= \frac{\pi^2 e(3e-2)}{6(e-1)^2} < 9.5.$$

For some positive integer $q$ and $R > 0$, we define $\mathcal{S}^1(q, R)$ as the set of continuous functions $f$ on $[0, 1)$ of the form

$$f(x) = \sum_{i=1}^{q+1}(\alpha_i x + \beta_i)\mathbb{1}_{[a_{i-1}, a_i)}(x)$$

with $0 = a_0 < a_1 < \cdots < a_{q+1} = 1$, $(\beta_1, \ldots, \beta_{q+1})' \in \mathbb{R}^{q+1}$ and $(\alpha_1, \ldots, \alpha_{q+1})' \in \mathbb{R}^{q+1}$, such that

$$\frac{1}{q}\sum_{i=1}^{q}|\alpha_{i+1} - \alpha_i| \leq R.$$

The following result holds.

COROLLARY 2. *Assume that $n \geq 9$. Let $K > 1$, $\kappa \in ]0,1[$, $\kappa' > 0$ and $p$ such that*

(5.6) $$p \leq (\kappa n - 2) \wedge (n - 9).$$

*Let $f \in \mathcal{S}^1(q, R)$ with $q \in \{1, \ldots, p\}$ and $R \leq \sigma e^{\kappa' n/q}$. If $\mu$ is defined by (2.1) then there exists a constant $C$ depending on $K$ and $\kappa, \kappa'$ only such that*

$$\mathbb{E}[\|\mu - \hat{\mu}_{\hat{m}}\|^2] \leq Cq\sigma^2\left[1 + \log\left(1 \vee \frac{nR^2}{q\sigma^2}\right)\right].$$

We postpone the proof of this result to Section 10.3.

5.5. *Estimating a signal.* We deal with the collection introduced in Section 2.4 and to each $m = (r, k_1, \ldots, k_d) \in \mathcal{M}$, associate the weight $L_m = (r+1)^d k_1 \cdots k_d$. With such a choice of weights, one can show that $\Sigma' \leq (e/(e-1))^{2(d+1)}$. For $\underline{\alpha} = (\alpha_1, \ldots, \alpha_d)$ and $\underline{R} = (R_1, \ldots, R_d)$ in $]0, +\infty[^d$, we denote by $\mathcal{H}(\underline{\alpha}, \underline{R})$ the space of $(\underline{\alpha}, \underline{R})$-Hölderian functions on $[0,1)^d$, which is the set of functions $f:[0,1)^d \to \mathbb{R}$ such that for any $i = 1, \ldots, d$ and $t_1, \ldots, t_d, z_i \in [0,1)$

$$\left|\frac{\partial^{r_i}}{\partial t_i^{r_i}} f(t_1, \ldots, t_i, \ldots, t_n) - \frac{\partial^{r_i}}{\partial t_i^{r_i}} f(t_1, \ldots, z_i, \ldots, t_n)\right| \leq R_i |t_i - z_i|^{\beta_i},$$

where $r_i + \beta_i = \alpha_i$, with $r_i \in \mathbb{N}$ and $0 < \beta_i \leq 1$.

In the sequel, we set $\|x\|_n^2 = \|x\|^2/n$ for $x \in \mathbb{R}^n$. By applying our procedure with the above weights and some $K > 1$, we obtain the following result.

COROLLARY 3. *Assume $n \geq 14$. Let $\underline{\alpha}$ and $\underline{R}$ fulfill the two conditions*

$$n^\alpha R_i^{2\alpha+d} \geq R^d \sigma^{2\alpha} \quad \text{and} \quad n^\alpha R_i^d \geq 2^\alpha R^d(r+1)^{d\alpha}, \quad \text{for } i = 1, \ldots, d,$$

*where*

$$r = \sup_{i=1,\ldots,d} r_i, \qquad \alpha = \left(\frac{1}{d}\sum_{i=1}^d \frac{1}{\alpha_i}\right)^{-1} \quad \text{and} \quad R = (R_1^{\alpha/\alpha_1}, \ldots, R_d^{\alpha/\alpha_d})^{1/d}.$$

*Then, there exists some constant $C$ depending on $r$ and $d$ only, such that for any $\mu$ given by (2.1) with $f \in \mathcal{H}(\underline{\alpha}, \underline{R})$,*

$$\mathbb{E}[\|\mu - \hat{\mu}\|_n^2] \leq C\left[\left(\frac{R^{d/\alpha}\sigma^2}{n}\right)^{2\alpha/(2\alpha+d)} \vee \left(\frac{R^2}{n^{2\alpha/d}}\right)\right].$$

The rate $n^{-2\alpha/(2\alpha+d)}$ is known to be minimax for density estimation in $\mathcal{H}(\underline{\alpha}, \underline{R})$ [see Ibragimov and Khas'minskii (1981)].

**6. Simulation study.** In order to evaluate the practical performance of our criterion, we carry out two simulation studies. In the first study, we consider the problem of detecting nonzero mean components. For the sake of comparison, we also include the performances of AIC, BIC and AMDL whose theoretical properties have been studied in Section 3. In the second study, we consider the variable selection problem and compare our procedure with adaptive Lasso recently proposed by Zou (2006). From a theoretical point of view, this last method cannot be compared with ours because its properties are shown assuming that the error variance is known. Nevertheless, this method gives good results in practice and the comparison with ours may be of interest. The calculations are made with R (www.r-project.org) and are available on request. We also mention that a simulation study has been carried out for the problem of multiple change-points detection (see Section 2.3). The results are available in Baraud, Giraud and Huet (2007).

6.1. *Computation of the penalties.* The calculation of the penalties we propose requires that of the EDkhi function or at least an upper bound for it. For $0 < q \leq 1$, the value $\mathsf{EDkhi}(D, N, q)$ is obtained by numerically solving for $x$ the equation

$$q = \mathbb{P}\left(F_{D+2,N} \geq \frac{x}{D+2}\right) - \frac{x}{D}\mathbb{P}\left(F_{D,N+2} \geq \frac{N+2}{ND}x\right),$$

where $F_{D,N}$ denotes a Fisher random variables with $D$ and $N$ degrees of freedom (see Lemma 6). However, this value of $x$ cannot be determined accurately enough when $q$ is too small. Rather, when $q < e^{-500}$ and $D \geq 2$, we bound the value of $\mathsf{EDkhi}(D, N, q)$ from above by solving for $x$ the equation

$$\frac{q}{2B(1+D/2, N/2)} = \frac{2 + NDx^{-1}}{N(N+2)}\left(\frac{N}{N+x}\right)^{N/2}\left(\frac{x}{N+x}\right)^{D/2},$$

where $B(p, q)$ stands for the beta function. This upper bound follows from formula (9.6), Lemma 6.

6.2. *Detecting nonzero mean components.*

*Description of the procedure.* We implement the procedure as described in Sections 2.1 and 5.2. More precisely, we select the set $\{(1), \ldots, (\hat{D})\}$ where $\hat{D}$ minimizes among $D$ in $\{1, \ldots, p\}$ the quantity defined at equation (5.3). In the case of our procedure, the penalty function $\overline{\mathrm{pen}}$ depends on a parameter $K$, and is equal to

$$\overline{\mathrm{pen}}_K(D) = K\frac{n-D}{n-D-1}\mathsf{EDkhi}\left[D+1, n-D-1, \left\{(D+1)^2\binom{n}{D}\right\}^{-1}\right].$$

We consider the three values $\{1; 1.1; 1.2\}$ for the parameter $K$ and denote $\hat{D}$ by $\hat{D}_K$, thus emphasizing the dependency on $K$. Even though the theory does not cover the case $K = 1$, it is worth studying the behavior of the procedure for this critical value. For the AIC, BIC and AMDL criteria, the penalty functions are respectively equal to

$$\overline{\text{pen}}_{\text{AIC}}(D) = (n - D)\left[\exp\left(\frac{2D}{n}\right) - 1\right],$$

$$\overline{\text{pen}}_{\text{BIC}}(D) = (n - D)\left[\exp\left(\frac{D\log(n)}{n}\right) - 1\right],$$

$$\overline{\text{pen}}_{\text{AMDL}}(D) = (n - D)\left[\exp\left(\frac{3D\log(n)}{n}\right) - 1\right].$$

We denote by $\hat{D}_{\text{AIC}}$, $\hat{D}_{\text{BIC}}$ and $\hat{D}_{\text{AMDL}}$ the corresponding values of $\hat{D}$.

*Simulation scheme.* For $\theta = (n, p, k, s) \in \mathbb{N} \times \{(p, k) \in \mathbb{N}^2 | k \leq p\} \times \mathbb{R}$, we denote by $\mathbb{P}_\theta$ the distribution of a Gaussian vector $Y$ in $\mathbb{R}^n$ whose components are independent with common variance 1 and mean $\mu_i = s$, if $i \leq k$ and $\mu_i = 0$ otherwise. Neither $s$ nor $k$ are known but we shall assume the upper bound $p$ on $k$ known:

$$\Theta = \{(2^j, p, k, s), j \in \{5, 9, 11, 13\}, p = \lfloor n/\log(n) \rfloor, k \in I_p, s \in \{3, 4, 5\}\},$$

where

$$I_p = \{2^{j'}, j' = 0, \ldots, \lfloor \log_2(p) \rfloor\} \cup \{0, p\}.$$

For each $\theta \in \Theta$, we evaluate the performance of each criterion as follows. On the basis of the 1000 simulations of $Y$ of law $\mathbb{P}_\theta$ we estimate the risk $R(\theta) = \mathbb{E}_\theta[\|\mu - \hat{\mu}_{\hat{m}}\|^2]$. Then, if $k$ is positive, we calculate the risk ratio $r(\theta) = R(\theta)/O(\theta)$, where $O(\theta)$ is the infimum of the risks over all $m \in \mathcal{M}$. More precisely,

$$O(\theta) = \inf_{m \in \mathcal{M}} \mathbb{E}_\theta[\|\mu - \hat{\mu}_m\|^2] = \inf_{D=0,\ldots,p}[s^2(k - D)\mathbf{I}_{D \leq k} + D].$$

It turns out that, in our simulation study, $O(\theta) = k$ for all $n$ and $s$.

*Results.* When $k = 0$, that is when the mean of $Y$ is 0, the results for AIC, BIC and AMDL criteria are given in Table 1. The theoretical results given in Section 3.2 and 3.3.2 are confirmed by the simulation study: when the complexity of the model collection $a$ equals $\log(n)$, AMDL satisfies the assumption of Theorem 1 and therefore the risk remains bounded, while the AIC and BIC criteria lead to an over-fitting (see Proposition 2). In all simulated samples, the BIC criterion selects a positive $\hat{D}$ and the AIC criterion chooses $\hat{D}$ equal to the largest possible dimension $p$. Our procedure, whose results are given in Table 2, performs similarly as AMDL. Since larger penalties tend to advantage small dimensional model, our procedure performs all the better that $K$ is large. AMDL overpenalizes models with positive dimension even more that $n$ is large, and then performs all the better.

TABLE 1
*Case $k = 0$. AIC, BIC and AMDL criteria: estimated risk R and percentage of the number of simulations for which $\widehat{D}$ is positive*

|       | AIC  |                          | BIC  |                          | AMDL |                            |
|-------|------|--------------------------|------|--------------------------|------|----------------------------|
| $n$   | $R$  | $\widehat{D}_{\text{AIC}} > 0$ | $R$  | $\widehat{D}_{\text{BIC}} > 0$ | $R$  | $\widehat{D}_{\text{AMDL}} > 0$ |
| 32    | 24   | 100%                     | 23   | 99%                      | 0.65 | 6.2%                       |
| 512   | 296  | 100%                     | 79   | 100%                     | 0.05 | 0.3%                       |
| 2048  | 1055 | 100%                     | 139  | 100%                     | 0.02 | 0.1%                       |
| 8192  | 3830 | 100%                     | 276  | 100%                     | 0.09 | 0.3%                       |

When $k$ is positive, Table 3 gives, for each $n$, the maximum of the risk ratios over $k$ and $s$. Note that the largest values of the risk ratios are achieved for the AMDL criterion. Besides, the AMDL risk ratio is maximum for large values of $k$. This is due to the fact that the quantity $3\log(n)$ involved in the AMDL penalty tends to penalize too severely models with large dimensions. Even in the favorable situation where the signal to noise ratio is large, AMDL criterion is unable to estimate $k$ when $k$ and $n$ are both too large. For example, Table 4 presents the values of the risk ratios when $k = n/16$ and $s = 5$, for several values $n$. Except in the situation where $n = 32$ and $k = 2$, the mean of the selected $\widehat{D}_{\text{AMDL}}$'s is small although the true $k$ is large. This overpenalization phenomenon is illustrated by Figure 1 which compares the AMDL penalty function with ours for $K = 1.1$. Let us now turn to the case where $k$ is small. The results for $k = 1$ are presented in Table 5. When $n = 32$, the methods are approximately equivalent whatever the value of $K$.

Finally, let us discuss the choice of $K$. When $k$ is large, the risk ratios do not vary with $K$ (see Table 4). Nevertheless, as illustrated by Table 5, $K$

TABLE 2
*Case $k = 0$. Estimated risk R and percentage of the number of simulations for which $\widehat{D}$ is positive using our penalty $\overline{\text{pen}}_K$*

|       | $K = 1$ |                  | $K = 1.1$ |                  | $K = 1.2$ |                  |
|-------|---------|------------------|-----------|------------------|-----------|------------------|
| $n$   | $R$     | $\widehat{D}_K > 0$ | $R$       | $\widehat{D}_K > 0$ | $R$       | $\widehat{D}_K > 0$ |
| 32    | 0.67    | 6.4%             | 0.40      | 3.7%             | 0.25      | 2.2%             |
| 512   | 0.98    | 5.7%             | 0.33      | 1.9%             | 0.07      | 0.4%             |
| 2048  | 1.00    | 5.1%             | 0.48      | 2.3%             | 0.09      | 0.4%             |
| 8192  | 0.96    | 4.2%             | 0.31      | 1.2%             | 0.14      | 0.5%             |

TABLE 3
*For each n, maximum of the estimated risk ratios $r_{\max}$ over the values of $(k,s)$ for $k > 0$. $\bar{k}$ and $\bar{s}$ are the values of $k$ and $s$ where the maxima are reached*

| | Our criterion with | | | | | | | | | AMDL | | |
|---|---|---|---|---|---|---|---|---|---|---|---|---|
| | $K = 1$ | | | $K = 1.1$ | | | $K = 1.2$ | | | | | |
| $n$ | $r_{\max}$ | $\bar{k}$ | $\bar{s}$ | $r_{\max}$ | $\bar{k}$ | $\bar{s}$ | $r_{\max}$ | $\bar{k}$ | $\bar{s}$ | $r_{\max}$ | $\bar{k}$ | $\bar{s}$ |
| 32 | 14.6 | 9 | 4 | 15.2 | 9 | 4 | 15.4 | 9 | 4 | 23.2 | 9 | 5 |
| 512 | 11.5 | 82 | 4 | 15.2 | 82 | 4 | 15.9 | 82 | 4 | 25.0 | 82 | 5 |
| 2048 | 10.7 | 1 | 4 | 15.5 | 268 | 4 | 16.0 | 256 | 4 | 25.0 | 256 | 5 |
| 8192 | 12.7 | 1 | 4 | 13.9 | 909 | 4 | 16.0 | 909 | 4 | 25.0 | 512 | 5 |

must stay close to 1 in order to avoid overpenalization. We suggest taking $K = 1.1$.

TABLE 4
*Case $k = n/16$ and $s = 5$. Estimated risk ratio $r$ and mean of the $\widehat{D}$'s*

| | | Our criterion with | | | | | | | AMDL | |
|---|---|---|---|---|---|---|---|---|---|---|
| | | $K = 1$ | | $K = 1.1$ | | $K = 1.2$ | | | | |
| $n$ | $k$ | $r$ | $\widehat{D}$ | $r$ | $\widehat{D}$ | $r$ | $\widehat{D}$ | | $r$ | $\widehat{D}$ |
| 32 | 2 | 3.43 | 2.04 | 3.89 | 1.94 | 4.49 | 1.85 | | 3.39 | 1.90 |
| 512 | 32 | 1.96 | 33.2 | 1.93 | 32.6 | 1.94 | 32.1 | | 23.5 | 2.12 |
| 2048 | 128 | 1.89 | 131 | 1.89 | 130 | 1.91 | 128 | | 25 | 0.52 |
| 8192 | 512 | 1.91 | 532 | 1.89 | 523 | 1.89 | 515 | | 25 | 0.22 |

TABLE 5
*Case $k = 1$ and $s = 5$. For each n, estimated risk ratio followed by the percentages of simulations for which $\widehat{D}$ is equal to $0, 1$ and larger than $1$*

| | Our criterion with | | | | | | | | | | | | |
|---|---|---|---|---|---|---|---|---|---|---|---|---|---|
| | $K = 1$ | | | | $K = 1.1$ | | | | $K = 1.2$ | | | | |
| | | Histogram | | | | Histogram | | | | Histogram | | | |
| $n$ | $R$ | $= 0$ | $= 1$ | $\geq 2$ | $R$ | $= 0$ | $= 1$ | $\geq 2$ | $R$ | $= 0$ | $= 1$ | $\geq 2$ | |
| 32 | 3.6 | 7.3 | 84.8 | 7.9 | 3.9 | 9.8 | 84.6 | 5.6 | 4.5 | 12.9 | 82.7 | 4.4 | |
| 512 | 5.4 | 14.6 | 80.4 | 5.0 | 6.1 | 20.3 | 77.8 | 1.9 | 7.2 | 26.0 | 73.0 | 1.0 | |
| 2048 | 7.1 | 21.8 | 74.9 | 3.3 | 8.2 | 28.6 | 70.1 | 1.3 | 9.6 | 35.4 | 64.1 | 0.5 | |
| 8192 | 9.1 | 29.5 | 67.7 | 2.8 | 10.4 | 37.4 | 61.6 | 1.0 | 12.2 | 45.9 | 53.9 | 0.2 | |

24     Y. BARAUD, C. GIRAUD AND S. HUET

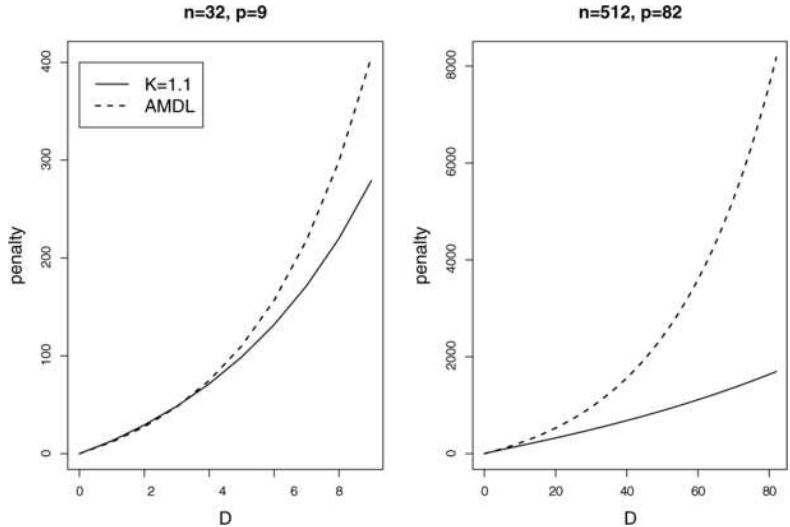

FIG. 1. *Comparison of the penalty functions* $\overline{\operatorname{pen}}_{\mathrm{AMDL}}(D)$ *and* $\overline{\operatorname{pen}}_K(D)$ *for* $K = 1.1$.

6.3. *Variable selection.* We present two simulation studies for illustrating the performances of our method for variable selection and compare them to the adaptive Lasso. The first simulation scheme was proposed by Zou ([2006](#)). The second one involves highly correlated covariates.

*Description of the procedure.* We consider the variable selection problem described in Section 2.2 and we implement the procedure considering the collection $\mathcal{M}$ for complete variable selection defined in Section 2.2.2 with maximal dimension $p$. We select the subset $\widehat{m}$ of $\{1, \ldots, N\}$ minimizing $\operatorname{Crit}_L(m)$ given at equation (1.4) with penalty function

$$\operatorname{pen}(m) = \overline{\operatorname{pen}}(|m|)$$
$$= K \frac{n - |m|}{n - |m| - 1} \mathsf{EDkhi}\bigg[|m| + 1, n - |m| - 1, \bigg\{p(|m| + 1)\binom{N}{|m|}\bigg\}^{-1}\bigg].$$

This choice for the penalty ensures a quasi oracle bound for the risk of $\widehat{m}$ [see inequality (5.5)].

*The adaptive Lasso procedure.* The adaptive Lasso procedure proposed by Zou starts with a preliminary estimator $\widehat{a}$ of $a$ as, for example, the ordinary least squares estimator when it exists. Then one computes the minimizer $\widehat{a}^w$ among those $a \in \mathbb{R}^N$ of the criterion

$$\operatorname{Crit}_{\mathrm{Lasso}}(a) = \bigg\|Y - \sum_{j=1}^N a_j x^{(j)}\bigg\|^2 + \lambda \sum_{j=1}^N \widehat{w}_j |a_j|,$$

where the weights $\widehat{w}_j = 1/|\widehat{a}|_j^\gamma$ for $j = 1, \ldots, N$. The smoothing parameters $\lambda$ and $\gamma$ are chosen by cross-validation. The set $\widehat{m}_{\text{Lasso}}$ is the set of indices $j$ such that $\widehat{a}_j^w$ is nonzero.

*Simulation scheme.* Let $M(n, N)$ be the set of matrices with $n$ rows and $N$ columns. For $\theta = (X, a, \sigma) \in M(n, N) \times \mathbb{R}^N \times \mathbb{R}^+$, we denote by $\mathbb{P}_\theta$ the distribution of a Gaussian vector $Y$ in $\mathbb{R}^n$ with mean $\mu = Xa$ and covariance $\sigma^2 I_n$. We consider two choices for the pair $(X, a)$. The first one is based on the *Model* 1 considered by Zou (2006) in its simulation study. More precisely, $N = 8$ and the rows of the matrix $X$ are $n$ i.i.d. Gaussian centered variables such that for all $1 \leq j < k \leq 8$ the correlation between $x^{(j)}$ and $x^{(k)}$ equals $0.5^{(k-j)}$. We did $S = 50$ simulations of the matrix $X$, denoted $\mathcal{X}^S = (X^s, s = 1, \ldots, S)$ and define

$$\Theta_1 = \{(X, a, \sigma), X \in \mathcal{X}^S, a = (3, 1.5, 0, 0, 2, 0, 0, 0)^T, \sigma \in \{1, 3\}\}.$$

The second one is constructed as follows. Let $x^{(1)}, x^{(2)}, x^{(3)}$ be three vectors of $\mathbb{R}^n$ defined by

$$x^{(1)} = (1, -1, 0, \ldots, 0)^T / \sqrt{2},$$
$$x^{(2)} = (-1, 1.001, 0, \ldots, 0)^T / \sqrt{1 + 1.001^2},$$
$$x^{(3)} = (1/\sqrt{2}, 1/\sqrt{2}, 1/n, \ldots, 1/n)^T / \sqrt{1 + (n-2)/n^2}$$

and for $4 \leq j \leq n$, let $x^{(j)}$ be the $j$th vector of the canonical basis of $\mathbb{R}^n$. We take $N = n$ and $\mu = (n, n, 0, \ldots, 0)^T$. Let $a \in \mathbb{R}^N$ satisfying $\mu = Xa$. Note that only the two first components of $a$ are nonzero. We thus define $\Theta_2 = \{(X, a, 1)\}$.

We choose $n = 20$ and for each $\theta \in \Theta_1 \cup \Theta_2$ we did 500 simulations of $Y$ with law $P_\theta$.

Our procedures were carried out considering all (nonvoid) subsets $m$ of $\{1, \ldots, N\}$ with cardinality not larger than $p = 8$. On the basis of the results obtained in the preceding section, we took $K = 1.1$.

For the adaptive Lasso procedure the parameters $\lambda$ and $\gamma$ are estimated using one-fold cross-validation as follows: when $\theta \in \Theta_1$, the values of $\lambda$ vary between 0 and 200 and following the recommendations given by Zou, $\gamma$ can take three values $(0.5, 1, 2)$. For $\theta \in \Theta_2$, $\lambda$ varies between 0 and 40, and $\gamma$ takes the values $(0.5, 1, 1.5)$; the value $\gamma = 2$ leading to numerical instability in the LARS algorithm.

We evaluate the performances of each procedure by estimating the risk ratio

$$r(\theta) = \frac{\mathbb{E}_\theta[\|\mu - \widehat{\mu}_{\widehat{m}}\|^2]}{\inf_{m \in \mathcal{M}} \mathbb{E}_\theta[\|\mu - \widehat{\mu}_m\|^2]},$$

the expectation of $|\widehat{m}|$, and calculating the frequencies of choosing and containing the true model $m_0$.

TABLE 6
*Case $\theta \in \Theta_1$. Risk ratio $r$, expectation of $|\widehat{m}|$ and percentages of the number of times $\widehat{m}$ equals or contains the true model ($m_0 = \{1, 2, 5\}$). These quantities are averaged over the $S$ design matrices $X$ in $\Theta_1$*

|           | $\sigma = 1$ |                       |                       |                           | $\sigma = 3$ |                       |                       |                           |
|-----------|------|----------------------|-------------------|---------------------------|------|----------------------|-------------------|---------------------------|
|           | $r$  | $\mathbb{E}(|\widehat{m}|)$ | $\widehat{m} = m_0$ | $\widehat{m} \supseteq m_0$ | $r$  | $\mathbb{E}(|\widehat{m}|)$ | $\widehat{m} = m_0$ | $\widehat{m} \supseteq m_0$ |
| $K = 1.1$ | 1.64 | 3.44                 | 67%               | 98.3%                     | 2.89 | 2.23                 | 12.4%             | 20.2%                     |
| A. Lasso  | 1.92 | 3.73                 | 62%               | 98.9%                     | 2.58 | 3.74                 | 13.7%             | 49.3%                     |

*Results.* When $\theta \in \Theta_1$, the methods give similar results. Looking carefully at the results shown in Table 6, we remark that the adaptive Lasso method selects more variables than ours. It gives results slightly better when $\sigma = 3$, the risk ratio being smaller and the frequency of containing the true model being greater. But, when $\sigma = 1$, using the adaptive Lasso method leads to increase the risk ratio and to wrongly detect a larger number of variables.

In case $\theta \in \Theta_2$, the adaptive Lasso procedure does not work while our procedure gives satisfactory results (see Table 7). The good behavior of our method in this case illustrates the strength of Theorem 2 whose results do not depend on the correlation of the explanatory variables.

Finally, let us emphasize that these methods are not comparable either from a theoretical point of view nor from a practical one. In our method the penalty function is free from $\sigma$, while in the adaptive Lasso method the theoretical results are given for known $\sigma$ and the penalty function depends on $\sigma$ through the parameter $\lambda$. All the difficulty of our method lies in the complexity of the collection $\mathcal{M}$, making impossible to consider in practice models with a large number of variables.

**7. Estimating the pair $(\mu, \sigma^2)$.** Unlike the previous sections which focused on the estimation of $\mu$, we consider here the problem of estimating the pair $\theta = (\mu, \sigma^2)$. All along, we shall assume that $\mathcal{M}$ is finite and consider

TABLE 7
*Case $\theta \in \Theta_2$ with $\sigma = 1$. Risk ratio $r$, expectation of $|\widehat{m}|$ and percentages of the number of times $\widehat{m}$ equals or contains the true model ($m_0 = \{1, 2\}$)*

|           | $r$  | $\mathbb{E}(|\widehat{m}|)$ | $\widehat{m} = m_0$ | $\widehat{m} \supseteq m_0$ |
|-----------|------|----------------------------|---------------------|-----------------------------|
| $K = 1.1$ | 2.35 | 2.28                       | 80.2%               | 96.6%                       |
| A. Lasso  | 26.5 | 10.2                       | 0.4%                | 40%                         |

the Kullback loss defined between $P_{\mu,\sigma^2}$ and $P_{\nu,\tau^2}$ by

$$\mathcal{K}(P_{\mu,\sigma^2}, P_{\nu,\tau^2}) = \frac{n}{2}\left[\log\left(\frac{\tau^2}{\sigma^2}\right) + \frac{\sigma^2}{\tau^2} - 1 + \frac{\|\mu - \nu\|^2}{n\tau^2}\right].$$

Given some finite collection of models $\mathcal{S} = \{S_m, m \in \mathcal{M}\}$ we associate to each $m \in \mathcal{M}$ the estimator $\hat{\theta}_m$ of $\theta$ defined by

$$\hat{\theta}_m = (\hat{\mu}_m, \hat{\sigma}_m^2) = \left(\Pi_m Y, \frac{\|Y - \Pi_m Y\|^2}{N_m}\right).$$

For a given $m$, the risk of $\hat{\theta}_m$ can be evaluated as follows.

PROPOSITION 5. *Let $\theta_m = (\mu_m, \sigma_m^2)$ where $\sigma_m^2 = \sigma^2 + \|\mu - \mu_m\|^2/n$ and $\mu_m = \Pi_m \mu$. Then,*

$$(7.1) \quad \inf_{\nu \in S_m, \tau^2 > 0} \mathcal{K}(P_\theta, P_{\nu,\tau^2}) = \mathcal{K}(P_\theta, P_{\theta_m}) = \frac{n}{2}\log\left(1 + \frac{\|\mu - \mu_m\|^2}{n\sigma^2}\right)$$

*and provided that $N_m > 2$,*

$$(7.2) \quad \mathbb{E}_\theta[\mathcal{K}(P_\theta, P_{\hat{\theta}_m})] \leq \mathcal{K}(P_\theta, P_{\theta_m}) + \frac{n}{2}\left[\frac{D_m + 2}{N_m - 2} - \log\left(1 - \frac{D_m}{n}\right)\right],$$

$$(7.3) \quad \mathbb{E}_\theta[\mathcal{K}(P_\theta, P_{\hat{\theta}_m})] \geq \mathcal{K}(P_\theta, P_{\theta_m}) \vee \left(\frac{N_m \wedge D_m}{2}\right).$$

*In particular, if $D_m \leq N_m$ and $N_m > 2$, then*

$$(7.4) \quad \begin{aligned} \mathcal{K}(P_\theta, P_{\theta_m}) \vee \frac{D_m}{2} &\leq \mathbb{E}[\mathcal{K}(P_\theta, P_{\hat{\theta}_m})] \\ &\leq \mathcal{K}(P_\theta, P_{\theta_m}) + 4(D_m + 2). \end{aligned}$$

As expected, this proposition shows that the Kullback risk of the estimator $\hat{\theta}_m$ is of order of a bias term, namely $\mathcal{K}(P_\theta, P_{\theta_m})$, plus some variance term which is proportional to $D_m$, at least when $D_m \leq (n/2) \wedge (n-3)$. We refer to Baraud, Giraud and Huet (2007) for the proof of these bounds.

Let us now introduce a definition.

DEFINITION 4. Let $F_{D,N}$ be a Fisher random variable with $D \geq 1$ and $N \geq 3$ degrees of freedom. For $x \geq 0$, we set

$$\mathsf{Fish}[D, N, x] = \frac{\mathbb{E}[(F_{D,N} - x)_+]}{\mathbb{E}(F_{D,N})} \leq 1.$$

For $0 < q \leq 1$ we define $\mathsf{EFish}[D, N, q]$ as the solution to the equation $\mathsf{Fish}[D, N, \mathsf{EFish}[D, N, q]] = q$.

We shall use the convention $\mathsf{EFish}[D, N, q] = 0$ for $q > 1$. Note that the restriction $N \geq 3$ is necessary to ensure that $\mathbb{E}(F_{D,N}) < \infty$.

Given some penalty $\text{pen}^*$ from $\mathcal{M}$ into $\mathbb{R}_+$, we shall deal with the penalized criterion

$$\text{Crit}'_K(m) = \frac{n}{2} \log\left(\frac{\|Y - \Pi_m Y\|^2}{N_m}\right) + \frac{1}{2} \text{pen}^*(m) \tag{7.5}$$

for which our results will take a more simple form than with criteria (1.4) and (1.5). In the sequel, we define

$$\tilde{\theta} = \hat{\theta}_{\hat{m}} \qquad \text{where } \hat{m} = \arg\min_{m \in \mathcal{M}} \text{Crit}'_K(m).$$

THEOREM 3. *Let $\mathcal{S} = \{S_m, m \in \mathcal{M}\}$, $\alpha = \min\{N_m/n | m \in \mathcal{M}\}$ and $K_1, K_2$ be two numbers satisfying $K_2 \geq K_1 > 1$. If $D_m \leq n - 5$ for all $m \in \mathcal{M}$, then the estimator $\tilde{\theta}$ satisfies*

$$\mathbb{E}[\mathcal{K}(P_\theta, P_{\tilde{\theta}})] \tag{7.6}$$
$$\leq \frac{K_1}{K_1 - 1}\left\{\inf_{m \in \mathcal{M}}\left[\mathbb{E}[\mathcal{K}(P_\theta, P_{\hat{\theta}_m})] + \frac{9}{4}(\text{pen}^*(m) \vee D_m)\right] + \Sigma_1 + \Sigma_2\right\},$$

*where*

$$\Sigma_1 = 2.5 e^{1/(K_2^2 \alpha)} n e^{-n/(4K_2^2)} |\mathcal{M}|^{4/(\alpha n)}, \qquad \Sigma_2 = \frac{5K_1}{4} \sum_{m \in \mathcal{M}}(D_m + 1)\Lambda_m$$

*and*

$$\Lambda_m = \mathsf{Fish}\left[D_m + 1, N_m - 1, \frac{N_m - 1}{K_1 N_m} \frac{K_2 D_m + (K_2 - 1)\text{pen}^*(m)}{K_2(D_m + 1)}\right].$$

*In particular, let $\mathcal{L} = \{L_m, m \in \mathcal{M}\}$ be a sequence of nonnegative weights. If for all $m \in \mathcal{M}$, $\text{pen}^*(m) = \text{pen}^{\mathcal{K}}_{K_1, K_2, \mathcal{L}}(m)$ with*

$$\text{pen}^{\mathcal{K}}_{K_1, K_2, \mathcal{L}}(m) = \frac{K_2}{K_2 - 1}\left[\frac{K_1(D_m + 1)N_m}{N_m - 1} \right. \tag{7.7}$$
$$\left. \times \mathsf{EFish}(D_m + 1, N_m - 1, e^{-L_m}) - D_m\right]_+,$$

*then the estimator $\tilde{\theta}$ satisfies (7.6) with $\Sigma_2 \leq 1.25 K_1 \sum_{m \in \mathcal{M}}(D_m + 1)e^{-L_m}$.*

This result is an analogue of Theorem 2 for the Kullbach risk. The expression of $\Sigma$ is akin to that of Theorem 2 apart from the additional term of order $ne^{-n/(4K_2^2)}|\mathcal{M}|^{4/(\alpha n)}$. In most of the applications, the cardinalities $|\mathcal{M}|$ of the collections are not larger than $e^{Cn}$ for some universal constant $C$, so that this additional term usually remains under control.

An upper bound for the penalty $\text{pen}^{\mathcal{K}}_{K_1, K_2, \mathcal{L}}$ is given in the following proposition, the proof of which is delayed to Section 10.2.

PROPOSITION 6. *Let $m \in \mathcal{M}$, with $D_m \geq 1$ and $N_m \geq 9$. We set $D = D_m + 1$, $N = N_m - 1$ and*

$$\Delta' = \frac{L_m + \log 5 + 1/(N-2)}{1 - 5/(N-2)}.$$

*Then, we have the following upper bound on the penalty $\mathrm{pen}^{\mathcal{K}}_{K_1,K_2,\mathcal{L}}$:*

(7.8) $\quad \mathrm{pen}^{\mathcal{K}}_{K_1,K_2,\mathcal{L}}(m) \leq \dfrac{K_1 K_2}{K_2 - 1} \dfrac{N+1}{N-2} \left[ 1 + e^{2\Delta'/N} \sqrt{\left(1 + \dfrac{2D}{N}\right) \dfrac{2\Delta'}{D}} \right]^2 D.$

## 8. Proofs of Theorems 2 and 3.

8.1. *Proof of Theorem 2.* We write henceforth $\varepsilon_m = \Pi_m \varepsilon$ and $\mu_m = \Pi_m \mu$. Expanding the squared Euclidean loss of the selected estimator $\hat{\mu}_{\hat{m}}$ gives

$$\begin{aligned}
\|\mu - \hat{\mu}_{\hat{m}}\|^2 &= \|\mu - \mu_{\hat{m}}\|^2 + \sigma^2 \|\varepsilon_{\hat{m}}\|^2 \\
&= \|\mu\|^2 - \|\mu_{\hat{m}}\|^2 + \sigma^2 \|\varepsilon_{\hat{m}}\|^2 \\
&= \|\mu\|^2 - \|\hat{\mu}_{\hat{m}}\|^2 + 2\sigma^2 \|\varepsilon_{\hat{m}}\|^2 + 2\sigma \langle \mu_{\hat{m}}, \varepsilon \rangle.
\end{aligned}$$

Let $m^*$ be an arbitrary index in $\mathcal{M}$. It follows from the definition of $\hat{m}$ that it also minimizes over $\mathcal{M}$ the criterion $\mathrm{Crit}(m) = -\|\hat{\mu}_m\|^2 + \mathrm{pen}(m)\hat{\sigma}_m^2$ and we derive

(8.1)
$$\begin{aligned}
\|\mu - \hat{\mu}_{\hat{m}}\|^2 &\leq \|\mu\|^2 - \|\hat{\mu}_{m^*}\|^2 + \mathrm{pen}(m^*)\hat{\sigma}_{m^*}^2 \\
&\quad - \mathrm{pen}(\hat{m})\hat{\sigma}_{\hat{m}}^2 + 2\sigma^2 \|\varepsilon_{\hat{m}}\|^2 + 2\sigma \langle \mu_{\hat{m}}, \varepsilon \rangle \\
&\leq \|\mu - \mu_{m^*}\|^2 - \sigma^2 \|\varepsilon_{m^*}\|^2 - 2\sigma \langle \mu_{m^*}, \varepsilon \rangle + \mathrm{pen}(m^*)\hat{\sigma}_{m^*}^2 \\
&\quad - \mathrm{pen}(\hat{m})\hat{\sigma}_{\hat{m}}^2 + 2\sigma^2 \|\varepsilon_{\hat{m}}\|^2 + 2\sigma \langle \mu_{\hat{m}}, \varepsilon \rangle \\
&\leq \|\mu - \mu_{m^*}\|^2 + R(m^*) - \mathrm{pen}(\hat{m})\hat{\sigma}_{\hat{m}}^2 + 2\sigma^2 \|\varepsilon_{\hat{m}}\|^2 \\
&\quad - 2\sigma \langle \mu - \mu_{\hat{m}}, \varepsilon \rangle,
\end{aligned}$$

where for all $m \in \mathcal{M}$,

$$R(m) = -\sigma^2 \|\varepsilon_m\|^2 + 2\sigma \langle \mu - \mu_m, \varepsilon \rangle + \mathrm{pen}(m)\hat{\sigma}_m^2.$$

For each $m$, we bound $\langle \mu - \mu_m, \varepsilon \rangle$ from above by using the inequality

(8.2) $\quad -2\sigma \langle \mu - \mu_m, \varepsilon \rangle \leq \dfrac{1}{K} \|\mu - \mu_m\|^2 + K\sigma^2 \langle u_m, \varepsilon \rangle^2,$

where $u_m = \mu - \mu_m / \|\mu - \mu_m\|$ when $\|\mu - \mu_m\| \neq 0$ and $u_m$ is any unit vector orthogonal to $S_m$ otherwise. Note that in any case, $\langle u_m, \varepsilon \rangle$ is a standard Gaussian random variable independent of $\|\varepsilon_m\|^2$. For each $m$, let $F_m$ be the linear space both orthogonal to $S_m$ and $u_m$. We bound $\hat{\sigma}_m^2$ from below by the following inequality:

$$N_m \frac{\hat{\sigma}_m^2}{\sigma^2} \geq \|\Pi_{F_m}\varepsilon\|^2, \tag{8.3}$$

where $\Pi_{F_m}$ denotes the orthogonal projector onto $F_m$.

By using (8.2), (8.3) and the fact that $2 - 1/K \leq K$, inequality (8.1) leads to

$$\frac{K-1}{K}\|\mu - \hat{\mu}_{\hat{m}}\|^2$$
$$\leq \|\mu - \mu_{m^*}\|^2 + R(m^*)$$
$$- \mathrm{pen}(\hat{m})\hat{\sigma}_{\hat{m}}^2 + (2 - 1/K)\sigma^2\|\varepsilon_{\hat{m}}\|^2 + K\sigma^2\langle u_{\hat{m}}, \varepsilon\rangle^2$$
$$\leq \|\mu - \mu_{m^*}\|^2 + R(m^*)$$
$$+ \sum_{m \in \mathcal{M}} [K\sigma^2\|\varepsilon_m\|^2 + K\sigma^2\langle u_m, \varepsilon\rangle^2 - \mathrm{pen}(m)\hat{\sigma}_m^2]\mathbf{1}_{\hat{m}=m}$$
$$\leq \|\mu - \mu_{m^*}\|^2 + R(m^*) + \sigma^2 \sum_{m \in \mathcal{M}} \left[KU_m - \mathrm{pen}(m)\frac{V_m}{N_m}\right]\mathbf{1}_{\hat{m}=m}, \tag{8.4}$$

where $U_m = \|\varepsilon_m\|^2 + \langle u_m, \varepsilon\rangle^2$ and $V_m = \|\Pi_{F_m}\varepsilon\|^2$. Note that $U_m$ and $V_m$ are independent and distributed as $\chi^2$ random variables with respective parameters $D_m + 1$ and $N_m - 1$.

8.1.1. *Case $c = 0$.* We start with the (simple) case $c = 0$. Then, by taking the expectation on both sides of (8.4), we get

$$\frac{K-1}{K}\mathbb{E}[\|\mu - \hat{\mu}_{\hat{m}}\|^2]$$
$$\leq \|\mu - \mu_{m^*}\|^2 + \mathbb{E}(R(m^*))$$
$$+ K\sigma^2 \sum_{m \in \mathcal{M}} \mathbb{E}\left(\left[U_m - \frac{(N_m - 1)\mathrm{pen}(m)}{KN_m} \times \frac{V_m}{N_m - 1}\right]_+\right)$$
$$\leq \|\mu - \mu_{m^*}\|^2 + \mathbb{E}(R(m^*))$$
$$+ K\sigma^2 \sum_{m \in \mathcal{M}} (D_m + 1)\mathsf{Dkhi}\left(D_m + 1, N_m - 1, \frac{(N_m - 1)\mathrm{pen}(m)}{KN_m}\right).$$

To conclude, we note that

$$\mathbb{E}(R(m^*)) = -\sigma^2 D_{m^*} + \mathrm{pen}(m^*)\left(\sigma^2 + \frac{\|\mu - \mu_{m^*}\|^2}{N_{m^*}}\right)$$

and $m^*$ is arbitrary among $\mathcal{M}$.

8.1.2. *Case $c > 0$.* We now turn to the case $c > 0$. We set $\bar{V}_m = V_m/N_m$ and $a_m = \mathbb{E}(\bar{V}_m)$. Analyzing the cases $\bar{V}_m \leq a_m$ and $\bar{V}_m > a_m$ apart gives

$$KU_m - \text{pen}(m)\bar{V}_m = [KU_m - (\text{pen}(m) + c - c)\bar{V}_m]\mathbf{1}_{\bar{V}_m \leq a_m}$$
$$+ [KU_m - (\text{pen}(m) + c - c)\bar{V}_m]\mathbf{1}_{\bar{V}_m > a_m}$$
$$\leq ca_m + [KU_m - (\text{pen}(m) + c)\bar{V}_m]_+ \mathbf{1}_{\bar{V}_m \leq a_m}$$
$$+ [KU_m - (\text{pen}(m) + c)a_m]_+ \mathbf{1}_{\bar{V}_m > a_m}$$
$$\leq c + [KU_m - (\text{pen}(m) + c)\bar{V}_m]_+$$
$$+ [KU_m - (\text{pen}(m) + c)\mathbb{E}(\bar{V}_m)]_+,$$

where we used for the final steps $a_m = \mathbb{E}(\bar{V}_m) \leq 1$. Going back to the bound (8.4), we obtain in the case $c > 0$

(8.5)
$$\frac{K-1}{K}\|\mu - \hat{\mu}_{\hat{m}}\|^2 \leq \|\mu - \mu_{m^*}\|^2 + R(m^*) + c\sigma^2$$
$$+ \sigma^2 \sum_{m \in \mathcal{M}} [KU_m - (\text{pen}(m) + c)\bar{V}_m]_+$$
$$+ \sigma^2 \sum_{m \in \mathcal{M}} [KU_m - (\text{pen}(m) + c)\mathbb{E}(\bar{V}_m)]_+.$$

Now, the independence of $U_m$ and $\bar{V}_m$ together with Jensen's inequality ensures that

$$\mathbb{E}([KU_m - (\text{pen}(m) + c)\mathbb{E}(\bar{V}_m)]_+) \leq \mathbb{E}([KU_m - (\text{pen}(m) + c)\bar{V}_m]_+),$$

so taking expectation in (8.5) gives

$$\frac{K-1}{K}\mathbb{E}[\|\mu - \hat{\mu}_{\hat{m}}\|^2]$$
$$\leq \|\mu - \mu_{m^*}\|^2 + \mathbb{E}(R(m^*)) + c\sigma^2$$
$$+ 2K\sigma^2 \sum_{m \in \mathcal{M}} \mathbb{E}\left(\left[KU_m - (\text{pen}(m) + c)\frac{V_m}{N_m}\right]_+\right)$$
$$\leq \|\mu - \mu_{m^*}\|^2 + \mathbb{E}(R(m^*)) + c\sigma^2$$
$$+ 2K\sigma^2 \sum_{m \in \mathcal{M}} (D_m + 1)\text{Dkhi}\left(D_m + 1, N_m - 1, \frac{(N_m - 1)(\text{pen}(m) + c)}{KN_m}\right).$$

To conclude, we follow the same lines as in the case $c = 0$.

8.2. *Proof of Theorem 3.* Let $m$ be arbitrary in $\mathcal{M}$. In the sequel we write $\mathcal{K}(m)$ for the Kullback divergence $\mathcal{K}(P_{\mu,\sigma^2}, P_{\hat{\mu}_m,\hat{\sigma}_m^2})$, namely

(8.6) $$\mathcal{K}(m) = \frac{n}{2}\log(\hat{\sigma}_m^2) + \frac{\|\mu - \hat{\mu}_m\|^2 + n\sigma^2}{2\hat{\sigma}_m^2} - \frac{n}{2}(\log\sigma^2 + 1).$$

We also set $\phi(x) = \log(x) + x^{-1} - 1 \geq 0$ for all $x \geq 0$, $\delta = 1/K_2$, and for each $m$ we define the random variable $\xi_m$ as the number $\langle u_m, \varepsilon \rangle$ with $u_m = \mu - \mu_m / \|\mu - \mu_m\|$ when $\|\mu - \mu_m\| \neq 0$ and $u_m$ is any unit vector orthogonal to $S_m$ otherwise.

We split the proof of Theorem 3 into four lemmas.

LEMMA 1. *The index $\hat{m}$ satisfies*

$$(8.7) \quad \frac{K_1 - 1}{K_1} \mathcal{K}(\hat{m}) \leq \mathcal{K}(m) + \frac{1 - \delta}{2} \mathrm{pen}^*(m) + R_1(m) + F(\hat{m}) + R_2(m, \hat{m})$$

*where, for all $m, m' \in \mathcal{M}$,*

$$R_2(m, m') = \frac{n(1-\delta) - \|\varepsilon\|^2}{2} \left( \frac{\sigma^2}{\hat{\sigma}_{m'}^2} - \frac{\sigma^2}{\hat{\sigma}_m^2} \right),$$

$$R_1(m) = \frac{D_m}{2} - \frac{\sigma^2 \|\varepsilon_m\|^2}{\hat{\sigma}_m^2} + \frac{\sigma \langle \mu - \mu_m, \varepsilon \rangle}{\hat{\sigma}_m^2} - \frac{\delta n}{2} \phi\left( \frac{\hat{\sigma}_m^2}{\sigma^2} \right),$$

$$F(m) = -\frac{D_m}{2} + \left( 1 - \frac{1}{2K_1} \right) \frac{\sigma^2 \|\varepsilon_m\|^2}{\hat{\sigma}_m^2} + \frac{K_1}{2} \frac{\sigma^2 \xi_m^2}{\hat{\sigma}_m^2} \mathbb{1}_{\{\xi_m < 0\}}$$
$$- \frac{1 - \delta}{2} \mathrm{pen}^*(m).$$

PROOF. We have

$$\mathcal{K}(\hat{m}) = \mathcal{K}(m) + \frac{n}{2} \log \frac{\hat{\sigma}_{\hat{m}}^2}{\hat{\sigma}_m^2} + \frac{\|\mu - \hat{\mu}_{\hat{m}}\|^2 + n\sigma^2}{2\hat{\sigma}_{\hat{m}}^2} - \frac{\|\mu - \hat{\mu}_m\|^2 + n\sigma^2}{2\hat{\sigma}_m^2}$$

$$= \mathcal{K}(m) + \frac{n}{2} \log \frac{\hat{\sigma}_{\hat{m}}^2}{\hat{\sigma}_m^2} + \frac{\|\mu - \hat{\mu}_{\hat{m}}\|^2 + n\sigma^2 - \|Y - Y_{\hat{m}}\|^2}{2\hat{\sigma}_{\hat{m}}^2} - \frac{D_{\hat{m}}}{2}$$

$$+ \frac{D_m}{2} - \frac{\|\mu - \hat{\mu}_m\|^2 + n\sigma^2 - \|Y - Y_m\|^2}{2\hat{\sigma}_m^2}$$

$$= \mathcal{K}(m) + \frac{n}{2} \log \frac{\hat{\sigma}_{\hat{m}}^2}{\hat{\sigma}_m^2} + \frac{2\|\varepsilon_{\hat{m}}\|^2 + n - \|\varepsilon\|^2}{2\hat{\sigma}_{\hat{m}}^2 / \sigma^2} - \frac{\sigma \langle \mu - \mu_{\hat{m}}, \varepsilon \rangle}{\hat{\sigma}_{\hat{m}}^2}$$

$$- \frac{D_{\hat{m}}}{2} + \frac{D_m}{2} - \frac{2\|\varepsilon_m\|^2 + n - \|\varepsilon\|^2}{2\hat{\sigma}_m^2 / \sigma^2} + \frac{\sigma \langle \mu - \mu_m, \varepsilon \rangle}{\hat{\sigma}_m^2}.$$

With $\xi_m$ defined before the lemma, we get

$$\mathcal{K}(\hat{m}) \leq \mathcal{K}(m) + \frac{n}{2} \log \frac{\hat{\sigma}_{\hat{m}}^2}{\hat{\sigma}_m^2} + \frac{2\|\varepsilon_{\hat{m}}\|^2 + n - \|\varepsilon\|^2}{2\hat{\sigma}_{\hat{m}}^2 / \sigma^2} + \frac{\|\mu - \mu_{\hat{m}}\|^2}{2K_1 \hat{\sigma}_{\hat{m}}^2} - \frac{D_{\hat{m}}}{2}$$

$$+ \frac{K_1 \mathbb{1}_{\{\xi_{\hat{m}} < 0\}} \xi_{\hat{m}}^2}{2\hat{\sigma}_{\hat{m}}^2 / \sigma^2} + \frac{D_m}{2} - \frac{2\|\varepsilon_m\|^2 + n - \|\varepsilon\|^2}{2\hat{\sigma}_m^2 / \sigma^2} + \frac{\sigma \langle \mu - \mu_m, \varepsilon \rangle}{\hat{\sigma}_m^2}.$$

In view of (8.6), since $\delta = 1/K_2 \leq 1/K_1 < 1$, we have

$$\frac{\|\mu - \mu_{\hat{m}}\|^2}{2K_1 \hat{\sigma}_{\hat{m}}^2} = \frac{\mathcal{K}(\hat{m})}{K_1} - \frac{\sigma^2 \|\varepsilon_{\hat{m}}\|^2}{2K_1 \hat{\sigma}_{\hat{m}}^2} - \frac{n}{2K_1} \phi\left(\frac{\hat{\sigma}_{\hat{m}}^2}{\sigma^2}\right)$$

$$\leq \frac{\mathcal{K}(\hat{m})}{K_1} - \frac{\sigma^2 \|\varepsilon_{\hat{m}}\|^2}{2K_1 \hat{\sigma}_{\hat{m}}^2} - \frac{\delta n}{2} \phi\left(\frac{\hat{\sigma}_{\hat{m}}^2}{\sigma^2}\right),$$

and thus,

$$\frac{K_1 - 1}{K_1} \mathcal{K}(\hat{m}) \leq \mathcal{K}(m) + \frac{n}{2} \log \frac{\hat{\sigma}_{\hat{m}}^2}{\hat{\sigma}_m^2} + \left(1 - \frac{1}{2K_1}\right) \frac{\sigma^2 \|\varepsilon_{\hat{m}}\|^2}{\hat{\sigma}_{\hat{m}}^2} - \frac{\delta n}{2} \phi\left(\frac{\hat{\sigma}_{\hat{m}}^2}{\sigma^2}\right)$$

$$+ \frac{K_1 \mathbf{1}_{\{\xi_{\hat{m}} < 0\}}}{2} \frac{\sigma^2 \xi_{\hat{m}}^2}{\hat{\sigma}_{\hat{m}}^2} - \frac{D_{\hat{m}}}{2} + \frac{D_m}{2} + \frac{n - \|\varepsilon\|^2}{2}\left(\frac{\sigma^2}{\hat{\sigma}_{\hat{m}}^2} - \frac{\sigma^2}{\hat{\sigma}_m^2}\right)$$

$$- \frac{\sigma^2 \|\varepsilon_m\|^2}{\hat{\sigma}_m^2} + \frac{\sigma \langle \mu - \mu_m, \varepsilon \rangle}{\hat{\sigma}_m^2}$$

$$\leq \mathcal{K}(m) + (1 - \delta) \frac{n}{2} \log \frac{\hat{\sigma}_{\hat{m}}^2}{\hat{\sigma}_m^2} + \left(1 - \frac{1}{2K_1}\right) \frac{\sigma^2 \|\varepsilon_{\hat{m}}\|^2}{\hat{\sigma}_{\hat{m}}^2}$$

$$+ \frac{K_1 \mathbf{1}_{\{\xi_{\hat{m}} < 0\}}}{2} \frac{\sigma^2 \xi_{\hat{m}}^2}{\hat{\sigma}_{\hat{m}}^2} - \frac{D_{\hat{m}}}{2} + R_2(m, \hat{m}) + R_1(m).$$

Finally, we get the result since $\hat{m}$ satisfies by definition $n \log(\hat{\sigma}_{\hat{m}}^2 / \hat{\sigma}_m^2) \leq \text{pen}^*(m) - \text{pen}^*(\hat{m})$. $\square$

LEMMA 2. *For all $m \in \mathcal{M}$, we have $\mathbb{E}(R_1(m)) \leq D_m/2$.*

PROOF. Since $\phi$ is nonnegative, we have

(8.8)
$$R_1(m) = \frac{D_m}{2} - \frac{\sigma^2 \|\varepsilon_m\|^2}{\hat{\sigma}_m^2} + \frac{\sigma \langle \mu - \mu_m, \varepsilon \rangle}{\hat{\sigma}_m^2} - \frac{\delta n}{2} \phi\left(\frac{\hat{\sigma}_m^2}{\sigma^2}\right)$$

$$\leq \frac{D_m}{2} + \frac{\sigma \langle \mu - \mu_m, \varepsilon \rangle}{\hat{\sigma}_m^2}.$$

Since $\varepsilon$ and $-\varepsilon$ have the same distribution, note that

$$2\mathbb{E}\left(\frac{\langle \mu - \mu_m, \varepsilon \rangle}{\hat{\sigma}_m^2}\right)$$

$$= (n - D_m) \mathbb{E}\left(\frac{\langle \mu - \mu_m, \varepsilon \rangle}{\|\mu - \mu_m\|^2 + \|\varepsilon - \varepsilon_m\|^2 + 2\langle \mu - \mu_m, \varepsilon \rangle}\right)$$

$$+ (n - D_m) \mathbb{E}\left(\frac{-\langle \mu - \mu_m, \varepsilon \rangle}{\|\mu - \mu_m\|^2 + \|\varepsilon - \varepsilon_m\|^2 - 2\langle \mu - \mu_m, \varepsilon \rangle}\right)$$

$$= (n - D_m)\mathbb{E}\left(\frac{-4\langle \mu - \mu_m, \varepsilon\rangle^2}{(\|\mu - \mu_m\|^2 + \|\varepsilon - \varepsilon_m\|^2)^2 - 4\langle \mu - \mu_m, \varepsilon\rangle^2}\right)$$

$$\leq 0.$$

Consequently, the result follows by taking the expectation on both sides of (8.8). □

LEMMA 3. *Under the assumptions that for all $m \in \mathcal{M}$, $N_m \geq \alpha n \geq 5$, we have for all $m \in \mathcal{M}$*

$$\mathbb{E}[R_2(m, \hat{m})] \leq \tfrac{7}{4}\mathrm{pen}^*(m) + 2.5 n e^{-(\alpha n - 4)\delta^2/(4\alpha)}|\mathcal{M}|^{4/(\alpha n)}.$$

PROOF.  Note that $R_2(m, \hat{m}) \leq R_{2,1}(m, \hat{m}) + R_{2,2}(m, \hat{m})$, where

$$R_{2,1}(m, \hat{m}) = \frac{1}{2}(\|\varepsilon\|^2 - (1 - \delta)n)_+ \left(\frac{\sigma^2}{\hat{\sigma}_m^2} - \frac{\sigma^2}{\hat{\sigma}_{\hat{m}}^2}\right)$$

and

$$R_{2,2}(m, \hat{m}) = \frac{1}{2}((1 - \delta)n - \|\varepsilon\|^2)_+ \frac{\sigma^2}{\hat{\sigma}_{\hat{m}}^2}.$$

It remains to bound the expectation of these two terms.

It follows from the definition of $\hat{m}$ and the inequality $1 - e^{-u} \leq u$ which holds for all $u \geq 0$ that

$$\frac{\sigma^2}{\hat{\sigma}_m^2} - \frac{\sigma^2}{\hat{\sigma}_{\hat{m}}^2} = \frac{\sigma^2}{\hat{\sigma}_m^2}\left(1 - \frac{\hat{\sigma}_m^2}{\hat{\sigma}_{\hat{m}}^2}\right)$$

$$\leq \frac{\sigma^2}{\hat{\sigma}_m^2}(1 - e^{-\mathrm{pen}^*(m)/n})$$

$$\leq \frac{\mathrm{pen}^*(m)}{n}\frac{\sigma^2}{\hat{\sigma}_m^2}$$

and thus,

$$\mathbb{E}[R_{2,1}(m, \hat{m})] = \frac{1}{2}\mathbb{E}\left[(\|\varepsilon\|^2 - (1 - \delta)n)_+ \left(\frac{\sigma^2}{\hat{\sigma}_m^2} - \frac{\sigma^2}{\hat{\sigma}_{\hat{m}}^2}\right)\right]$$

$$\leq \frac{1}{2}\mathbb{E}\left([\|\varepsilon\|^2 - (1 - \delta)n]_+ \frac{\sigma^2}{\hat{\sigma}_m^2}\right)\frac{\mathrm{pen}^*(m)}{n}$$

$$\leq \frac{1}{2}\mathbb{E}([\|\varepsilon\|^2 - (1 - \delta)n]^2)^{1/2}\mathbb{E}\left(\frac{\sigma^4}{\hat{\sigma}_m^4}\right)^{1/2}\frac{\mathrm{pen}^*(m)}{n}$$

$$\leq \frac{\sqrt{\delta^2 + 2/n}}{2}\frac{N_m}{\sqrt{(N_m - 2)(N_m - 4)}}\mathrm{pen}^*(m)$$

$$\leq \tfrac{7}{4}\mathrm{pen}^*(m).$$

As to $\mathbb{E}[R_{2,2}(m,\hat{m})]$, we apply Hölder's inequality with $p = \lfloor \alpha n/4 \rfloor + 1$, $q = p/(p-1)$ and have

$$\mathbb{E}[R_{2,2}(m,\hat{m})] = \frac{1}{2}\mathbb{E}\left[(n(1-\delta) - \|\varepsilon\|^2)_+ \frac{\sigma^2}{\hat{\sigma}_{\hat{m}}^2}\right]$$

$$\leq \frac{n}{2}\mathbb{E}\left[\frac{\sigma^2}{\hat{\sigma}_{\hat{m}}^2}\mathbb{1}_{\|\varepsilon\|^2 \leq n(1-\delta)}\right]$$

$$\leq \frac{n}{2}[\mathbb{P}(\|\varepsilon\|^2 \leq n(1-\delta))]^{1/q}\mathbb{E}\left(\frac{\sigma^{2p}}{\hat{\sigma}_{\hat{m}}^{2p}}\right)^{1/p}$$

$$\leq \frac{n}{2}[\mathbb{P}(\|\varepsilon\|^2 \leq n(1-\delta))]^{1/q}\left[\sum_{m \in \mathcal{M}} \mathbb{E}\left(\frac{\sigma^{2p}}{\hat{\sigma}_m^{2p}}\right)\right]^{1/p},$$

and by using that $\mathbb{P}(\|\varepsilon\|^2 \leq n(1-\delta)) \leq \exp(-n\delta^2/4)$ [see Laurent and Massart (2000), Lemma 1] together with (9.2) (note that $N_{m'} > 2p$ for all $m' \in \mathcal{M}$)

$\mathbb{E}[R_{2,2}(m,\hat{m})]$

$$\leq \frac{n}{2}e^{-n\delta^2/(4q)}\left[\sum_{m \in \mathcal{M}} \frac{N_m^p}{(N_m - 2)(N_m - 4)\cdots(N_m - 2p)}\right]^{1/p}$$

$$\leq \frac{n}{2}e^{-n\delta^2/(4q)}\left[\sum_{m \in \mathcal{M}} \frac{N_m^p}{(N_m - 2)(N_m - 4)\cdots(N_m - 2p)}\right]^{1/p}$$

$$\leq 2.5n\sigma^2 e^{-n\delta^2/(4q)}|\mathcal{M}|^{1/p} \leq 2.5ne^{-(\alpha n - 4)\delta^2/(4\alpha)}|\mathcal{M}|^{4/(\alpha n)}. \quad \square$$

LEMMA 4. *Under the assumption that $N_m \geq 5$ for all $m \in \mathcal{M}$, we have*

(8.9) $$\mathbb{E}[F(\hat{m})] \leq \frac{5K_1}{4}\sum_{m \in \mathcal{M}}(D_m + 1)\mathsf{Fish}[D_m + 1, N_m - 1, q_m]$$

*with*

$$q_m = \frac{(N_m - 1)}{K_1(D_m + 1)N_m}\left[D_m + \frac{K_2 - 1}{K_2}\mathrm{pen}^*(m)\right].$$

PROOF. Since $\mathbb{E}[F(\hat{m})] \leq \sum_{m \in \mathcal{M}} \mathbb{E}[F(m)]$, it suffices to bound $\mathbb{E}[F(m)]$ from above for all $m$. As in the proof of Theorem 2, we introduce $U_m = \|\varepsilon_m\|^2 + \xi_m^2$ and $V_m = \|\Pi_{F_m}\varepsilon\|^2 \leq N_m\hat{\sigma}_m^2/\sigma^2$. Since $\delta = 1/K_2$, we get

$$F(m) = \left[\left(1 - \frac{1}{2K_1}\right)\|\varepsilon_m\|^2 + \frac{K_1}{2}\mathbb{1}_{\{\xi_m < 0\}}\xi_m^2\right]\frac{\sigma^2}{\hat{\sigma}_m^2} - \frac{1}{2}(D_m + (1-\delta)\mathrm{pen}^*(m))$$

$$\leq \frac{K_1}{2}\frac{N_m U_m}{V_m} - \frac{1}{2}\left(D_m + \frac{K_2-1}{K_2}\text{pen}^*(m)\right)$$

$$\leq \frac{K_1 N_m(D_m+1)}{2(N_m-1)}$$

$$\times \left[\frac{U_m(N_m+1)}{V_m(D_m-1)} - \frac{N_m-1}{K_1(D_m+1)N_m}\left(D_m + \frac{K_2-1}{K_2}\text{pen}^*(m)\right)\right].$$

Since $\frac{U_m(N_m+1)}{V_m(D_m-1)}$ is distributed as a Fisher random variable with $D_m+1$ and $N_m-1$ degrees of freedom, the result follows by taking the expectation on both sides and using $N_m \geq 5$. □

*End of the proof of Theorem 3.* By taking the expectation on both sides of (8.7) and using Lemmas 2, 3 and 4 (we recall that $\delta = 1/K_2$) we obtain

$$\frac{K_1-1}{K_1}\mathbb{E}[\mathcal{K}(\hat{m})]$$

$$\leq \mathbb{E}[\mathcal{K}(m)] + \frac{9}{4}\text{pen}^*(m) + \frac{D_m}{2} + \frac{5}{2}ne^{-(\alpha n-4)\delta^2/(4\alpha)}|\mathcal{M}|^{4/(\alpha n)}$$

$$+ \frac{5K_1}{4}\sum_{m'\in\mathcal{M}}(D_{m'}+1)\text{Fish}[D_{m'}+1, N_{m'}-1, q_{m'}],$$

which leads to (7.6) since $m$ is arbitrary in $\mathcal{M}$. Note that the latter series is not larger than $\sum_{m'\in\mathcal{M}}(D_{m'}+1)e^{-L_{m'}}$ for $\text{pen}^*(m') = \text{pen}^{\mathcal{K}}_{K_1,K_2,\mathcal{L}}(m)$ by definition of EFish.

**9. Some preliminary results.** The aim of this section is to establish some technical results we shall use hereafter. The proofs of these being elementary and mainly based on integration by parts, we omit them and rather refer the interested reader to the technical report Baraud, Giraud and Huet (2007). We start with some moment inequalities on the inverse of a $\chi^2$ random variable.

LEMMA 5. *Let $V$ be a $\chi^2$ random variable with $N > 2$ degrees of freedom and noncentrality parameter $a$. We have*

(9.1) $$\frac{1}{a+N-2} \leq \mathbb{E}\left(\frac{1}{V}\right) \leq \frac{N}{(N+a)(N-2)} \leq \frac{1}{N-2}.$$

*Let $p$ be some positive integer. If $N > 2p$, then*

(9.2) $$\mathbb{E}\left(\frac{1}{V^p}\right) \leq \frac{1}{(N-2)\cdots(N-2p)}.$$

*Besides, equality holds in (9.2) for $a = 0$.*

We recall that $\phi(t) = (t - 1 - \log(t))/2$ for all $t \geq 1$. For two positive integers $D$ and $N$, $F_{D,N}$ denotes a Fisher random variable with $D$ and $N$ degrees of freedom, and we set

$$B(N/2, D/2) = \int_0^1 t^{N/2}(1-t)^{D/2}\, dt, \tag{9.3}$$

$$\psi_{D,N}(t) = \phi(t) - \frac{D(t-1)^2}{4(D+N+2)} \qquad \text{for all } t \geq 1. \tag{9.4}$$

The following holds.

LEMMA 6. *Let $D$ and $N$ be two positive integers. For all $x \geq 0$,*

$$\mathsf{Dkhi}(D, N, x) = \mathbb{P}\left(F_{D+2,N} \geq \frac{x}{D+2}\right) \tag{9.5}$$
$$- \frac{x}{D}\mathbb{P}\left(F_{D,N+2} \geq \frac{(N+2)x}{DN}\right).$$

*If $D \geq 2$ and $x \geq D$, then*

$$\mathsf{Dkhi}(D, N, x)$$
$$\leq \frac{1}{B(N/2, 1+D/2)}\left(\frac{N}{N+x}\right)^{N/2}\left(\frac{x}{N+x}\right)^{D/2}\frac{2(2x+ND)}{N(N+2)x} \tag{9.6}$$

$$\leq \left(1 + \frac{2x}{ND}\right)\mathbb{P}\left(F_{D,N+2} \geq \frac{(N+2)x}{ND}\right) \tag{9.7}$$

$$\leq \left(1 + \frac{2x}{ND}\right)\exp\left[-D\psi_{D,N}\left(\frac{(N+2)x}{ND}\right)\right]. \tag{9.8}$$

The next lemma states similar bounds on $\mathsf{Fish}(D, N, x)$.

LEMMA 7. *Let $D$ and $N$ be integer fulfilling $D \geq 1$ and $N \geq 3$. Then, for any $x \geq 0$,*

$$\mathsf{Fish}(D, N, x) \tag{9.9}$$
$$= \mathbb{P}\left(F_{D+2,N-2} \geq \frac{(N-2)D}{(D+2)N}x\right) - x\frac{N-2}{N}\mathbb{P}(F_{D,N} \geq x),$$

*where $F_{D,N}$ is a Fisher random variable with $D$ and $N$ degrees of freedom. Moreover, when $x \geq \frac{N}{N-2}$ and $D \geq 2$, we have the upper bounds*

$$\mathsf{Fish}(D, N, x)$$
$$\leq \frac{2}{B(D/2, N/2)}\left(\frac{N}{N+Dx}\right)^{N/2}\left(\frac{Dx}{N+Dx}\right)^{D/2-1}\frac{2x+N}{N^2} \tag{9.10}$$

$$\leq \left(\frac{1+2x}{N}\right)\mathbb{P}(F_{D,N} \geq x). \tag{9.11}$$

**10. Proofs of propositions and corollaries.**

10.1. *Proof of Proposition 3.* Let $m \in \mathcal{M}$, $D \in \{0,\ldots,n-2\}$, $N = n - D$ and $\mathcal{M}_D = \{m \in \mathcal{M}, D_m = D\}$. For all $c \geq 0$, (4.3) implies that $\{m' \in \mathcal{M}_D | L_{m'} \leq c\}$ is finite and since the map $x \mapsto \mathsf{EDkhi}(D+1, N-1, x)$ is decreasing, so is
$$\{m' \in \mathcal{M}_D | \mathsf{EDkhi}(D+1, N-1, e^{-L_{m'}}) \leq c\}.$$
It follows from the definitions of $\mathrm{Crit}_L$ and $\mathrm{pen}_{K,\mathcal{L}}$ that for some nonnegative constant $c = c(Y, D, n, m)$,
$$\bar{\mathcal{M}}_D = \{m' \in \mathcal{M}_D | \mathrm{Crit}_L(m') \leq \mathrm{Crit}_L(m)\}$$
is a subset of $\{m' \in \mathcal{M}_D | \mathsf{EDkhi}(D+1, N-1, e^{-L_{m'}} \leq c)\}$ and is therefore also finite. We deduce that $\mathrm{Crit}_L$ is minimum for some element of the finite set $\bar{\mathcal{M}} = \bigcup_{D=0}^{n-2} \bar{\mathcal{M}}_D$, thus showing that $\hat{m}$ exists. The remaining part of the proposition follows by taking $c = 0$ in Theorem 2.

10.2. *Proofs of Propositions 4 and 6.* Let us start with the proof of Proposition 4. We set
$$b(\Delta, D, N) = \left[1 + e^{2\Delta/(N+2)} \sqrt{\left(1 + \frac{2D}{N+2}\right) \frac{2\Delta}{D}}\right]^2$$
and $x = Db(\Delta, D, N) \geq D$.

Since $\mathrm{pen}_{K,\mathcal{L}}(m) = \frac{K(N+1)}{N} \mathsf{EDkhi}(D, N, e^{-L_m})$, we obtain (4.5) by showing the inequality $\mathsf{EDkhi}(D, N, e^{-L_m}) \leq x$ or equivalently

(10.1) $\qquad\qquad\qquad \mathsf{Dkhi}(D, N, x) \leq e^{-L_m}.$

Let us now turn to the proof of (10.1). Since $D \geq 2$ and $x \geq D$, we can apply (9.7) and get
$$\mathsf{Dkhi}(D, N, x) \leq \left(\frac{1 + 2x}{ND}\right) \mathbb{P}\left(F_{D,N+2} \geq \frac{(N+2)x}{ND}\right)$$
$$\leq \left(1 + \frac{2b(\Delta, D, N)}{N}\right) \mathbb{P}(F_{D,N+2} \geq b(\Delta, D, N)).$$

The deviation inequality on Fisher random variables available in Baraud, Huet and Laurent (2003) (Lemma 1) gives with $F = F_{D,N+2}$

$\mathbb{P}(F_{D,N+2} \geq b(\Delta, D, N))$
$$\leq \mathbb{P}\left(F \geq 1 + 2\sqrt{\left(1 + \frac{D}{N+2}\right) \frac{\Delta}{D}} + e^{4\Delta/(N+2)}\left(1 + \frac{2D}{N+2}\right) \frac{2\Delta}{D}\right)$$
$$\leq \mathbb{P}\left(F \geq 1 + 2\sqrt{\left(1 + \frac{D}{N+2}\right) \frac{\Delta}{D}} + \left(1 + \frac{2D}{N+2}\right) \frac{N+2}{2D}[e^{4\Delta/(N+2)} - 1]\right)$$
$$\leq e^{-\Delta}$$

and hence,

$$\mathsf{Dkhi}(D, N, x) \leq \left(1 + \frac{2b(\Delta, D, N)}{N}\right) e^{-\Delta}.$$

By using $D \geq 2$ and $N \geq 6$, we crudely bound $b(\Delta, D, N)$ from above as follows:

$$b(\Delta, D, N) = \left[1 + e^{2\Delta/(N+2)}\sqrt{\left(1 + \frac{2D}{N+2}\right)\frac{2\Delta}{D}}\right]^2$$

$$\leq \left[1 + \sqrt{\Delta}\sqrt{\frac{3}{2}e^{4\Delta/N}}\right]^2$$

$$\leq (1 + \Delta)\left(1 + \frac{3}{2}e^{4\Delta/N}\right)$$

$$\leq \frac{5}{2}(1 + \Delta)e^{4\Delta/N}$$

and deduce

$$\mathsf{Dkhi}(D, N, x) \leq \left(1 + \frac{2b(\Delta, D, N)}{N}\right) e^{-\Delta}$$

$$\leq 5e^{4\Delta/N}\left(1 + \frac{1+\Delta}{N}\right)e^{-\Delta}$$

$$\leq 5e^{1/N}e^{-\Delta(1-5/N)}.$$

Since $\Delta(1 - 5/N) = L_m + \log 5 + 1/N$, inequality (10.1) follows, thus completing the proof of (4.5).

We turn to (4.6). When $D_m = 0$, we obtain (4.6) by showing (10.1) for

$$x = 3\left[1 + e^{2L_m/N}\sqrt{\left(1 + \frac{6}{N}\right)\frac{2L_m}{3}}\right]^2.$$

We deduce from (9.5) that $\mathsf{Dkhi}(1, N, x) \leq \mathbb{P}(F_{3,N} \geq x/3)$. Again, the deviation inequality on Fisher random variables gives, with $L = L_m$,

$\mathbb{P}(F_{3,N} \geq x/3)$

$$= \mathbb{P}\left(F_{3,N} \geq \left[1 + e^{2L/N}\sqrt{\left(1 + \frac{6}{N}\right)\frac{2L}{3}}\right]^2\right)$$

$$\leq \mathbb{P}\left(F_{3,N} \geq 1 + 2\sqrt{\left(1 + \frac{3}{N}\right)\frac{L}{3}} + e^{4L/N}\left(1 + \frac{6}{N}\right)\frac{2L}{3}\right)$$

$$\leq \mathbb{P}\left(F_{3,N} \geq 1 + 2\sqrt{\left(1 + \frac{3}{N}\right)\frac{L}{3}} + \left(1 + \frac{6}{N}\right)\frac{N}{6}[e^{4L/N} - 1]\right) \leq e^{-L},$$

leading to (10.1). The proof of Proposition 4 is complete.

Since the proof of Proposition 6 is similar, we only give the main steps. We set

$$b'(\Delta', D, N) = \left[1 + e^{2\Delta'/N}\sqrt{\left(1+\frac{2D}{N}\right)\frac{2\Delta'}{D}}\right]^2$$

and $x' = Nb'(\Delta', D, N)/(N-2) \geq N/(N-2)$. In view of (9.11) and Lemma 1 in Baraud, Huet and Laurent (2003), we have

$$\mathsf{Fish}(D, N, x') \leq \left(1 + \frac{2b'(\Delta', D, N)}{N-2}\right)\mathbb{P}(F_{D,N} \geq b'(\Delta', D, N))$$

$$\leq \left(1 + \frac{2b'(\Delta', D, N)}{N-2}\right)e^{-\Delta'}.$$

Furthermore, when $D \geq 2$ and $N \geq 8$,

$$b'(\Delta', D, N) \leq [1 + e^{2\Delta/N}\sqrt{3/2}\sqrt{\Delta'}]^2 \leq \tfrac{5}{2}(1+\Delta')e^{4\Delta'/N},$$

which enforces

$$\mathsf{Fish}(D, N, x') \leq 5e^{1/(N-2)}e^{-\Delta'(1-5/(N-2))} \leq e^{-L_m}.$$

As a consequence,

$$\mathrm{pen}^{\mathcal{K}}_{K_1, K_2, \mathcal{L}}(m) \leq \frac{K_1 K_2}{K_2 - 1}\frac{D(N+1)}{N}\mathsf{EFish}(D, N, e^{-L_m}) \leq \frac{K_1 K_2}{K_2 - 1}\frac{D(N+1)}{N}x'$$

$$\leq \frac{K_1 K_2}{K_2 - 1}\frac{N+1}{N-2}\left[1 + e^{2\Delta'/N}\sqrt{\left(1+\frac{2D}{N}\right)\frac{2\Delta'}{D}}\right]^2 D.$$

10.3. *Proof of Corollary 2.* We start with an approximation lemma.

LEMMA 8. *For all $f \in \mathcal{S}^1(q, R)$ and $j \geq 1$ such that $1 \leq q \leq 2^j - 1$, there exists $m \in \mathcal{M}_{j,q}$ and $g \in \mathcal{F}_m$ such that $\|f - g\|_\infty \leq Rq2^{-j}$.*

PROOF. For $j \geq 1$ and $a \in [0, 1]$, we define $a^{(j)} = \inf\{x \in \mathcal{D}_j : x \geq a\}$. For all $x \in [0, 1)$, one can write

$$f(x) = f(0) + \int_0^x \sum_{i=1}^{q+1} \alpha_i \mathbb{1}_{[a_{i-1}, a_i)}(t)\, dt.$$

We take for $x \in [0, 1)$,

$$g(x) = f(0) + \int_0^x \sum_{i=1}^{q+1} \alpha_i \mathbb{1}_{[a_{i-1}^{(j)}, a_i^{(j)})}(t)\, dt.$$

Since one may have $a_{i-1}^{(j)} = a_i^{(j)}$ for some indices $i$, the function $g$ belongs to some space $\mathcal{F}_{m'}$ with $m' \in \mathcal{M}_{j,q'}$ and $q' \leq q$. By taking (any) $m \in \mathcal{M}_{j,q}$ such that $m' \subset m$, one has $g \in \mathcal{F}_m$.

For each $i \in \{1, \ldots, q+1\}$, we either have $a_{i-1} \leq a_{i-1}^{(j)} < a_i \leq a_i^{(j)}$ or $a_{i-1} < a_i \leq a_{i-1}^{(j)} = a_i^{(j)}$. In any case, we have

$$\mathbb{1}_{[a_{i-1}, a_i[} - \mathbb{1}_{[a_{i-1}^{(j)}, a_i^{(j)}[} = \mathbb{1}_{[a_{i-1}, a_{i-1}^{(j)}[} - \mathbb{1}_{[a_i, a_i^{(j)}[}$$

and consequently, for all $x \in [0, 1)$,

$$f(x) - g(x) = \int_0^x \sum_{i=1}^{q+1} \alpha_i (\mathbb{1}_{[a_{i-1}, a_i[}(t) - \mathbb{1}_{[a_{i-1}^{(j)}, a_i^{(j)}[}(t)) \, dt$$

$$= \int_0^x \sum_{i=1}^{q+1} \alpha_i (\mathbb{1}_{[a_{i-1}, a_{i-1}^{(j)}[}(t) - \mathbb{1}_{[a_i, a_i^{(j)}[}(t)) \, dt$$

$$= \int_0^x \left[ \sum_{i=1}^q (\alpha_{i+1} - \alpha_i) \mathbb{1}_{[a_i, a_i^{(j)}[}(t) \right.$$

$$\left. + \alpha_1 \mathbb{1}_{[a_0, a_0^{(j)}[}(t) - \alpha_{q+1} \mathbb{1}_{[a_{q+1}, a_{q+1}^{(j)}[}(t) \right] dt.$$

Since $a_0 = a_0^{(j)} = 0$, $a_{q+1} = a_{q+1}^{(j)} = 1$, and $|a_i^{(j)} - a_i| \leq 2^{-j}$, we obtain for $f \in \mathcal{S}^1(q, R)$

$$|f(x) - g(x)| \leq \sum_{i=1}^q |\alpha_{i+1} - \alpha_i| 2^{-j} \leq Rq 2^{-j}.$$

□

We take $m \in \mathcal{M}_{j,q}$ as in the lemma above with $j$ such that

$$2^{j-1} \leq \max\left\{ q, \sqrt{\frac{nR^2 q}{\sigma^2}} \right\} \leq 2^j.$$

We deduce from Proposition 4 [inequality (4.5)] that when $p \leq (\kappa n - 2) \wedge (n - 9)$ and $R \leq \sigma e^{\kappa' n / q}$, we have

$$\mathrm{pen}_{K, \mathcal{L}}(m) \leq C(K, \kappa) q \left( \frac{2^j}{q} \right)^{36q/(1-\kappa)n} \log\left( \frac{e 2^j}{q} \right)$$

$$\leq C(K, \kappa, \kappa') q \left[ 1 + \log\left( 1 \vee \frac{nR^2}{q \sigma^2} \right) \right].$$

Besides,

$$\frac{\|\mu - \mu_m\|^2}{\sigma^2} \left( 1 + \frac{\mathrm{pen}_{K, \mathcal{L}}(m)}{N_m} \right) \leq \frac{nR^2 q^2 2^{-2j}}{\sigma^2} \left( 1 + \frac{\mathrm{pen}_{K, \mathcal{L}}(m)}{N_m} \right)$$

$$\leq q\left(1 + \frac{\text{pen}_{K,\mathcal{L}}(m)}{N_m}\right)$$

$$\leq C'(K, \kappa, \kappa')q\left[1 + \log\left(1 \vee \frac{nR^2}{q\sigma^2}\right)\right],$$

and the result follows from (4.2).

10.4. *Proof of Corollary 3.* We write

$$\eta = \left(\frac{R^{d/\alpha}\sigma^2}{n}\right)^{\alpha/(2\alpha+d)} \vee \left(\frac{R(r+1)^\alpha}{(n/2)^{\alpha/d}}\right),$$

and set $m = (r, k_1, \ldots, k_d)$ where

$$k_i = \left\lfloor \left(\frac{R_i}{\eta}\right)^{1/\alpha_i} \right\rfloor, \qquad i = 1, \ldots, d.$$

It follows from our choice of $\eta \geq \frac{(R(r+1)^\alpha)}{(n/2)^{\alpha/d}}$ and the assumption $n \geq 14$ that

$$(r+1)^d k_1 \cdots k_d \leq n/2 \leq n - 2.$$

Moreover, under the assumptions $n^\alpha R_i^{2\alpha+d} \geq R^d \sigma^{2\alpha}$ and $n^{\alpha/d} R_i \geq 2^{\alpha/d} R(r+1)^\alpha$ we have $k_i \geq 1$ for all $i$. Consequently, $m \in \mathcal{M}$.

From formula (4.25) in Barron, Birgé and Massart (1999) we know that there exist a constant $C = C(d, r)$ and a piecewise polynomial $P$ in $\mathcal{F}_m$ such that

$$\|f - P\|_\infty \leq C \sum_{i=1}^d R_i k_i^{-\alpha_i} \leq C'\eta.$$

Moreover, since the assumptions of Proposition 4 hold, we have

$$\text{pen}_{K,\mathcal{L}}(m) \leq C(K)L_m \leq C(K)(r+1)^d R^{2d/(2\alpha+d)}(n/\sigma^2)^{d/(2\alpha+d)},$$

where the second inequality follows from the fact that $\eta \geq (\frac{R^{d/\alpha}\sigma^2}{n})^{\alpha/(2\alpha+d)}$. It remains to apply Theorem 2 to obtain the result.

**Acknowledgment.** We are very grateful to Lucien Birgé for his useful comments on an early draft of the present paper.

Y. BARAUD  
C. GIRAUD  
LABORATOIRE J. A. DIEUDONNÉ  
UNIVERSITÉ DE NICE SOPHIA ANTIPOLIS  
PARC VALROSE  
06108 NICE CEDEX 02  
FRANCE  
E-MAIL: yannick.baraud@unice.fr  
       Christophe.Giraud@polytechnique.edu

S. HUET  
INRA MIAJ  
78352 JOUY EN JOSAS CEDEX  
FRANCE  
E-MAIL: sylvie.huet@jouy.inra.fr